\documentclass{amsart}
\usepackage{amscd}

\newtheorem{lemma}{Lemma}[section]
\newtheorem{theorem}[lemma]{Theorem}
\newtheorem{remark}[lemma]{Remark}

\newtheorem{coro}[lemma]{Corollary}
\newtheorem{definition}[lemma]{Definition}
\newtheorem{example}[lemma]{Example}

\title[Almost Periodic motions and Attractors of
  Navier-Stokes Flows]{Almost Periodic Solutions and
Global Attractors of Non-autonomous Navier-Stokes Equations}
\author{David ~Cheban}
\address[D. Cheban]{%
State University of Moldova\\ Departement of Mathematics and
Informatics\\ A. Mateevich Street 60\\ MD--2009 Chi\c{s}in\u{a}u,
Moldova}
\email[D. Cheban]{cheban@usm.md}
\author{Jinqiao ~Duan}
\address[J. Duan]{%
Departement of Applied Mathematics \\ Illinois Institute of
Technology\\ Chicago, IL 60616, USA\\and Department of Mathematics\\
University of Science and Technology of China\\
Hefei, Anhui 230026, China}
\email[J.~Duan]{duan@iit.edu}
\date{\today}

\begin{document}
\begin{abstract}
{\bf J. Dynamics and Diff. Eqns., in press, 2004. }

 { The article is devoted to the study of non-autonomous
Navier-Stokes equations. First, the authors have proved that such
systems admit compact global attractors. This problem is
formulated and solved in the terms of general non-autonomous
dynamical systems. Second, they have obtained   conditions of
convergence of non-autonomous Navier-Stokes equations. Third, a
criterion for the existence of almost periodic (quasi
periodic,almost automorphic, recurrent, pseudo recurrent)
solutions of non-autonomous Navier-Stokes equations is given.
Finally, the authors have derived a global averaging principle for
non-autonomous Navier-Stokes equations. }
\end{abstract}

\maketitle
\section{Introduction}
We consider the two-dimensional Navier-Stokes system
\begin{eqnarray}\label{eq1.1}
& & u'+q(t)\sum _{i=1}^{2}u_i\partial _i u =\nu \Delta u -\nabla p + \phi (t) 
\nonumber \\
& & div\ u =0, \ \ u|_{\partial D} =0 ,
\end{eqnarray}
where $ D $ is an open bounded set with boundary $ \partial D \in C^{2}$.
This equation can be written in the following form
\begin{equation}\label{eq1.2}
u'+Au +B(t)(u,u) = f(t)
\end{equation}
on the corresponding Sobolev's space $E$, where $-A$ is a Stokes operator, 
$B(t)$ is a
bilinear form satisfying the identity
\begin{equation}\label{eq1.3}
Re \langle B(t)(u,v),w \rangle = -Re \langle B(t)(u,w),u \rangle
\end{equation}
for all $t\in \mathbb R $ and $u,v,w \in E$, and $f$ is forcing term.

In the work \cite{CV},\cite{DF},\cite{Il96},\cite{Il98} there
is studied a non-stationary equation (\ref{eq1.2}), when
$f$ is a function of time $t\in \mathbb R$. It is shown that the equation 
with compact
$f$ (in particularly, almost periodic) admits a compact global attractor and 
also
for small nonlinear (bilinear) term it was proved the existence a unique 
almost
periodic (quasi periodic, periodic) solution of equation (\ref{eq1.2}) if the 
forcing
term $f$ is almost periodic (quasi periodic, periodic).

The aim of the present article is to study the
equation
(\ref{eq1.2}) in the case, when the the bilinear form $ B $,
and the function $ f $ are non-stationary. The conditions
under which a non-stationary equation of type (\ref{eq1.2}) admits a
compact global attractor are indicated.

The theorem of "partial" averaging on finite interval for ordinary
differential equations it was proved in thew work \cite{Hap}. The
works \cite{DF},\cite{Il96} and \cite{Il98} are devoted to
generalization of method of averaging for dissipative partial
differential equations. We prove the theorem of "partial"
averaging for non-autonomous Navier-Stokes equation (\ref{eq1.2})
(i.e. the bilinear form and forcing term are non-stationaries).

Our paper is organized as follows:

In Section 2 we introduce a class of non-autonomous Navier-Stokes
equations and establish its dissipativity (Theorem \ref{t2.2}).

In Section 3 we prove that non-autonomous Navier-Stokes equations
admit a compact global attractor (Theorem \ref{t3.6}).

Section 4 is devoted to study of the problem of existence of
almost periodic (quasi periodic, recurrent, pseudo recurrent)
solutions of non-autonomous Navier-Stokes equations (Corollaries
\ref{cor4.7}) and we give the conditions of convergence of this
equations (Theorem \ref{t4.4}).

In Section 5 we prove the uniform averaging principle for the
non-autonomous Navier-Stokes equations on the finite segment
(Theorem \ref{t5.3}).

Section 6 is devoted to prove the global averaging principle for
non-autonomous Navier-Stokes equations on the semi-axis (Theorems
\ref{t6.1},\ref{t6.3} and \ref{t6.4}).

\section{Non-autonomous Navier-Stokes equations.}

Some results from the theory of semigroups of linear operators
\cite{Hen} and PDEs \cite{Il96}, \cite{SY},\cite{Tem} are
collected below.

A closed operator $A$ with domain $D(A)$ that is dense in a Banach
space $X$ is called a sectorial operator if for some $a\in \mathbb
R $ and $\varphi \in (0,\frac{\pi}{2})$ the sector
\begin{equation}
S_{a,\varphi}:=\{ \lambda \in \mathbb C, \pi \ge \vert
arg(\lambda - a)\vert \ge \varphi \}
\end{equation}
is contained in the resolvent set and for $\lambda \in
S_{a,\varphi} $
\begin{equation}
\Vert (\lambda I - A)^{-1}\Vert _{X\to X} \le \frac{c}{\vert
\lambda - a\vert +1}.
\end{equation}

For a sectorial operator $A$ the analytic semigroup of linear
bounded operators in $X$ is defined and denoted by $e^{-At}, \
t\ge 0$.

Let $A$ be a sectorial operator with $Re \sigma (A)
>0.$ For $\alpha \in (0,1)$ we define fractional powers of $A$ as
follows: $$ A^{\alpha}:=(A^{-\alpha})^{-1}, \ \text{where} \
A^{-\alpha}:= \frac{1}{\Gamma (\alpha)} \int
_{0}^{\infty}t^{\alpha -1}e^{-A t}dt. $$

The corresponding domains $D(A^{\alpha})$ are Banach spaces with norm given 
by
$$
\vert \cdot \vert _{\alpha} :=\vert \cdot \vert _{D(A^{\alpha})}
=\vert A^{\alpha}\cdot\vert .
$$

\begin{theorem}
The following estimates are valid:
\begin{enumerate}
\item
\begin{equation}
\Vert e^{-A t}\Vert_{X\to X} \le Ce^{-a t}, \ \ t \ge 0,
\end{equation}
\item
\begin{equation}
\Vert A^{\alpha}e^{-A t}\Vert_{X\to X} \le C_{\alpha}t^{-\alpha}
e^{-a t}, \ \ t > 0.
\end{equation}
\end{enumerate}
\end{theorem}

Let $ \Omega $ be a compact metric space, $ \mathbb R = (-\infty,
+\infty), (\Omega,\mathbb R, \sigma)$ be a dynamical system on $ \Omega $,
$\mathcal E$ be a real or complex Hilbert space, $ L(\mathcal E)$ be the 
space
of all linear forms on $ \mathcal E $, $L^{2}(\mathcal E)$ be the space of 
all bilinear
continuous forms $B:\mathcal E \times \mathcal E \to \mathcal F$\
and $ C(\Omega,W) $ be a space
of all continuous functions $ f: \Omega \to W\ $ ( $W$ is some metric
space), endowed with the topology of uniform convergence. Let us
consider the equation
\begin{equation}\label{eq2.1}
 u'+ Au + B(\omega t)(u,u) = f(\omega t),
\end{equation}
$(\omega \in \Omega ) $ where $\omega t:=\sigma(t,\omega),
B \in C(\Omega ,L^2(\mathcal E ))$, $ f
\in C(\Omega ,\mathcal E )$ and $A$ is a linear operator.

Below we will use some notions, denotations and results from \cite{Il98}.
Let Hilbert spaces $E, F, X$ satisfy $E \subset F;
 \ E, F, X \subset \mathcal E, $ each embedding being dense and continuous.

We further suppose that the linear operator $A$ is densely defined in
$\mathcal E$ and such that the linear equation
\begin{equation}\label{eq2.2}
 u'+ A u =0
\end{equation}
generates the $c_0$-semigroup of linear bounded operators
$$
e^{-At}: \mathcal E \to \mathcal E, \  \varphi (t,x):=e^{-At}x,
$$
which for $t>0$ can be extended to the linear bounded operators from $F$
to $E$ satisfying the following estimates
\begin{equation}\label{eq2.3}
\Vert e^{-At}\Vert_{E\to E} \le K e^{-a t},
\end{equation}
\begin{equation}\label{eq2.4}
\Vert e^{-At}\Vert_{F\to E} \le K t^{-\alpha_1} e^{-a t},
\ \ 0\le \alpha _1 <1 ,
\end{equation}
\begin{equation}\label{eq2.5}
\Vert Ae^{-At}\Vert_{F\to E} \le K t^{-\alpha_2} e^{-a t},
\ \ 0\le \alpha _2 <2 .
\end{equation}

We also suppose that the following condition is satisfied
\begin{equation}\label{eq2.5*}
 Ae^{At}=e^{At}A ,
\end{equation}
in the sense of $ L (F,E):=\{A:F \to E \ | A \ \text{is linear and
bounded} \ \}$ equipped with the operational norm.

Bilinear form $B.$ Denote by $ L^{2} (E,F)$ the space of all
bilinear bounded form $B: E\times E \to F$ with the norm
\begin{eqnarray}\label{eq5.8*}
 &&\Vert B \Vert :=\sup \{ \vert B(u,v) \vert _{F} \ : \ \vert u \vert \le 1,
\vert v \vert \le 1 \} . \nonumber
\end{eqnarray}
Let $C(\Omega ,L^{2} (E,F)) $ be a space of all continuous
mappings $B : \Omega \to L^{2} (E,F)$ and
\begin{eqnarray}\label{eq5.9*}
&& C_ B :=\sup \{ \vert B(\omega)(u,v) \vert _{F} \ : \  \omega \in \Omega ,
 \ \vert u \vert \le 1, \vert v \vert \le 1 \} ,\nonumber
\end{eqnarray}
then the mapping $F: \Omega \times E \to F \ (F(\omega ,u):=B(\omega)(u,u)) $
satisfies the following inequality
\begin{equation}\label{eq5.10*}
\vert B(\omega)(u_1,u_1)-B(\omega)(u_2,u_2) \vert _{F} \le C_B (\vert 
u_1\vert_{E} +
\vert u_2\vert_{E})\vert u_1-u_2 \vert _{E}
\end{equation}
for all $u_1,u_2 \in E .$

From the inequality (\ref{eq5.10*}) follows that on every ball $B[0,R]:=
\{ u\in E \ : \ \vert u \vert \le R \}$  we have
\begin{equation}\label{eq5.11*}
\vert B(\omega)(u_1,u_1)-B(\omega)(u_2,u_2) \vert _{F} \le
2C_B R \vert u_1-u_2 \vert _{E}
\end{equation}
for all $u_1,u_2 \in E .$

\begin{remark}\label{r5.1*}
The space of all the bilinear form $C(\Omega ,L^{2} (E,F)) $ is a
Banach space with the norm $\Vert B \Vert :=C_B .$
\end{remark}

Function $f$. The external force $f: \Omega \to X $ is continuous, i.e.
$f\in C(\Omega ,X)$.

Operators $e^{-At} .$ The operators $e^{-At} \ (t>0)$ can be
extended to the linear bounded operators from $X$ to $E$
satisfying the estimates
\begin{equation}\label{eq5.12*}
\Vert e^{-At}\Vert _{X\to E} \le Kt^{-\beta _1}e^{-at}, \ 0\le \beta _1 <1,
\end{equation}
\begin{equation}\label{eq5.13*}
\Vert Ae^{-At}\Vert _{X\to E} \le Kt^{-\beta _2}e^{-at}, \ 0\le \beta _2 <2,
\end{equation}
and the equation (\ref{eq2.5*}), this time in the sense of
$L(X,E)$ .

We suppose that the following conditions are fulfilled:
\begin{enumerate}
\item
there exists $\alpha > 0 $ such that
\begin{equation}\label{eq2.6}
Re \langle Au,u \rangle \ge  \alpha \vert u \vert^2
\end{equation}
for all $ u \in E ,$
where $ \vert \cdot \vert $ is a norm in $ E $ ;
\item
\begin{equation}\label{eq2.7}
Re \langle B(\omega )(u,v),w \rangle = - Re \langle B(\omega )(u,w),v \rangle
\end{equation}
for every $ u,v,w \in E $ and $ \omega \in \Omega$ .
\end{enumerate}

\begin{remark}\label{rem2.1}
a. It follows from (\ref{eq2.7}) that
\begin{equation}\label{eq2.8}
Re \langle B(\omega )(u,v),v) \rangle = 0
\end{equation}
for every $ u,v \in E$ and $ \omega \in \Omega .$

b.
\begin{equation}\label{eq2.9}
\vert B(\omega)(u,v) \vert _{F}\le C_{B} \vert u \vert_{E} \vert v \vert_{E}
\end{equation}
for all $ u,v \in E $ and $ \omega \in \Omega$, where $ C_B= \sup \{ \vert
B(\omega)(u,v) \vert_{F} : \omega \in \Omega,\ u,v\in E,\
\vert u \vert_{E} \le 1,$ and $
\vert v \vert_{E} \le 1 \}$.
\end{remark}

The equation (\ref{eq2.1}) with conditions (\ref{eq2.6}) and
(\ref{eq2.7}) is called a non-autonomous Navier-Stokes equation.
We will consider the mild solutions of the equation (\ref{eq2.1}),
i.e. $ u \in C([0,T],E)$ and satisfy the following integral
equation
\begin{equation}\label{eq^{*}}
u(t)=e^{-At}x +\int _0^{t}e^{-A(t-s)}(-B(\omega s)(u(s),u(s))+f(\omega s))ds.
\end{equation}

\begin{theorem}\label{t*}
Let $x_0\in E$, $r>0 $ and the conditions
(\ref{eq2.3}),(\ref{eq2.4}) and (\ref{eq2.6}) are fulfilled, then
there exist positive numbers $\delta =\delta (x_0,r)$ and
$T=T(x_0,r) $ such that the equation (\ref{eq^{*}}) admits a
unique solution $\varphi (t,x,\omega)$ ($x\in B[x_0,\delta]=\{x\in
E \ | \ \vert x -x_0 \vert \le \delta \} $) defined on the
interval $[0,T]$ with the conditions: $\varphi (0,x,\omega)=x$,
$\vert \varphi (t,x,\omega)-x_0\vert \le r$ for all $t\in [0,T]$
and the mapping $ \varphi : [0,T]\times B[x_0,\delta] \times
\Omega \to E \ ( (t,x,\omega)\to \varphi (t,x,\omega))$ is
continuous.
\end{theorem}
\begin{proof}
Let $x_0 \in E, \ r>0, \ \delta >0$ and $T>0$. We consider a space
$C_{x_0,r,\delta, T}$ of all continuous functions $\psi :
[0,T]\times B[x_0,\delta] \times \Omega \to B[x_0,r]$ equipped with the 
distance
$$
d(\psi _1,\psi_2):=\sup \{ \vert \psi _1(t,x,\omega)-\psi (t,x,\omega)
\vert _{E} \ : \ 0\le t \le T, x\in B[x_0,\delta], \omega \in \Omega \}
$$
is a complete metric space.

We define the operator $\Phi $ acting onto $C_{x_0,r,\delta, T} $
by the equality $$ (\Phi \psi)(t,x,\omega)=e^{-At}x +\int
_0^{t}e^{-A(t-s)} (-B(\omega
s)(\psi(s,x,\omega),\psi(s,x,\omega))+f(\omega s))ds . $$ There
exist $\delta _1=\delta _1(x_0,r)>0$ and $T_1=T_1(x_0,r)>0$ such
that $ \Phi C_{x_0,r,\delta, T} \subseteq C_{x_0,r,\delta, T}$ for
all $\delta \in (0,\delta _1]$ and $T \in (0,T _1]$. In fact,
\begin{eqnarray}\label{eq^{**}}
& & \vert  (\Phi \psi)(t,x,\omega)-x_0\vert _{E} \le
\vert  e^{-At}x-x_0\vert _{E} + \nonumber \\
&& \vert \int _0^{t}e^{-A(t-s)}
B(\omega s)(\psi(s,x,\omega),\psi(s,x,\omega)))ds \vert _{E} +
\vert \int _0^{t}e^{-A(t-s)}f(\omega s)ds  \vert _{E} \le \nonumber \\
& & m(\delta,T) + \int _0^{t}Ke^{-a(t-s)}(t-s)^{-\alpha _1}
\vert \psi (s,x,\omega) \vert _{E}^{2}ds +  \nonumber \\
& & \int _0^{t}Ke^{-a(t-s)}(t-s)^{-\beta _1}
\Vert f \Vert ds\le m(\delta ,T) +
K(\vert x_0 \vert _{E} +r)^{2}\frac{T^{1-\alpha _1}}{1-\alpha _1} \nonumber 
\\
& & + K\Vert f \Vert \frac{T^{1-\beta _1}}{1-\beta _1}:=d_1(x_0,r,\delta ,T)
\to 0  \nonumber
\end{eqnarray}
as $\delta +T\to 0,$
where $m(\delta,T):=\sup \{ \vert e^{-tA}x -x_0\vert _{E} \ : \
t\in [0,T], x\in B[x_0,r] $) and $\Vert f \Vert := \sup \{ \vert f(\omega)
\vert_{X} \ : \omega \in \Omega \}.$ Thus there exist $\delta _1=\delta 
_1(x_0,r)>0$ and
$T_1=T_1(x_0,r)>0$ such that $d_1(x_0,r,\delta,T)\le r$ for all
$\delta \in (0,\delta _1]$ and $T\in (0,T_1]$.

Let now $\psi _1, \psi _2 \in C_{x_0,r,\delta, T}$, then
\begin{eqnarray}\label{eq^{***}}
& & \vert (\Phi \psi _1)(t,x,\omega))-(\Phi \psi _2)(t,x,\omega))\vert _{E}= 
\nonumber \\
& & \vert  \int _{0}^{t} [B(\omega s)(\psi _1(s,x,\omega),\psi 
_1(s,x,\omega))-
B(\omega s)(\psi _2(s,x,\omega),\psi _2(s,x,\omega))]\vert _{E} \le 
\nonumber\\
& & 2C_B(\vert x_0 \vert _{E} +r)Td(\psi _1,\psi _2) \nonumber
\end{eqnarray}
and, consequently, $d(\Phi \psi_1,\Phi \psi_2)\le L(x_0,r,T)d(\psi
_1,\psi _2),$ where $L(x_0,r,T)=2C_B(\vert x_0\vert +r)T \to 0$ as
$T\to 0.$ Thus there exists $T_2=T_2(x_0,r)>0$ such that
$L(x_0,r,T)<1 $ for all $T\in (0,T_2].$ Denote by $\delta
(x_0,r):=\delta _1(x_0,r)$ and $T(x_0,r):=min(T_1(x_0,r),
T_2(x_0,r))$, then the mapping $\Phi :C_{x_0,r,\delta, T} \to
C_{x_0,r,\delta, T} $ is a contraction and, consequently, there
exists a unique function $\varphi \in C_{x_0,r,\delta, T} $
satisfying the equation (\ref{eq^{*}}) on the interval $[0,T]$.
The theorem is proved.
\end{proof}

\begin{remark}\label{rr*}
The theorem \ref{t*} is true and for the equation
$$
u'+Au=\mathcal F (\omega t,u)
$$
if the continuous function
$\mathcal F : \Omega \times E \to F$ satisfies the following
conditions:
\begin{enumerate}
\item
$$
\sup \{\vert \mathcal F (\omega,0)\vert _{E} \ : \ \omega \in \Omega \}
< \infty
$$
($\Omega$ ,generally speaking, is not compact);
\item
$F$ is locally Lipschitz, i.e. for every $r>0$ there exists $L(r)>0$ such 
that
$$
\vert \mathcal F(\omega, u_1) -\mathcal F(\omega, u_2) \vert _{F}
\le L(r) \vert u_1-u_2 \vert _{E}
$$
for all $u_1,u_2 \in E$ with condition: $\vert u_i\vert _{E} \le r \ 
(i=1,2).$
\end{enumerate}
\end{remark}

\begin{theorem}\label{t**}
Let $\mathcal K$ be a family of solutions of equation
(\ref{eq^{*}}) satisfying the following condition: there exists a
positive constant $M$ such that $\vert x(t)\vert _{\mathcal D
(A)}\le M $ for all $t\in \mathbb R_+ \ (\vert x \vert _{\mathcal
D (A)}:= \vert Ax\vert _{E}).$ If there exists $\tilde{C_{B}}>0$
such that
$$
\vert B(\omega)(u,v)\vert _{F} \le \tilde{C_{B}}\vert u \vert _{\mathcal D 
(A)}
\vert v \vert _{\mathcal D (A)}
$$
for all $u,v \in \mathcal D (A) $ , then this family of functions
is uniform equicontinuous on $\mathbb R_{+}$, i.e. for every
$\varepsilon >0$ there exists $\delta (\varepsilon)>0$ such that
$\vert t_1-t_2\vert <\delta $ implies $\vert x(t_1)-x(t_2)\vert
<\varepsilon $ for all $t_1,t_2 \in \mathbb R_{+}$ and $x\in
\mathcal K$.
\end{theorem}
\begin{proof} Let $\psi \in \mathcal K$ and $x:=\psi (0)$,
then $\psi (t)=\varphi (t,x,\omega)$
for all $t\in \mathbb R_+$ and we have
\begin{eqnarray}\label{eq2.10*}
\vert \varphi (t,x,\omega) -x\vert _{E} \le \vert e^{-At}x -x \vert _{E} +
\vert \int _0 ^{t}e^{-A(t-\tau)}(-B(\omega s)(\varphi (s,x,\omega),
\varphi (s,x,\omega)) + \nonumber \\
f(\omega s) )ds\vert _{E} \le \int _0^{t}e^{-as}s^{-\alpha _1}
\vert x \vert _{\mathcal D (A)}ds +
\int _0^{t}e^{-as}(t-s)^{-\alpha _1}C_b
\vert \varphi (s,x,\omega) \vert^{2} _{\mathcal D (A)}ds + \nonumber \\
\int _0^{t}e^{-a(t-s)}(t-s)^{-\beta _1}\Vert f \Vert ds \le
\frac{t^{1-\alpha_1}}{1-\alpha _1}M +C_BM^{2}\frac{t^{1-\alpha_1}}{1-\alpha 
_1} +
\Vert f \Vert \frac{t^{1-\beta_1}}{1-\beta _1}.
\end{eqnarray}
From (\ref{eq2.10*}) we obtain
$$
\sup \{ \vert \varphi (t,x,\omega)-x\vert_{E} \ : \ \vert x \vert _{\mathcal 
D (A)}, \
\omega \in \Omega \} \to 0
$$
as $t\to 0$ and, consequently,
\begin{eqnarray}\label{eq2.10**}
&& \vert \varphi (t_2,x,\omega) - \varphi (t_1,x,\omega)\vert _{E} =
\vert \varphi (t_2-t_1,\varphi (t_1,x,\omega), \omega t_1) -
\varphi (t_1,x,\omega)\vert _{E} \le \nonumber \\
&& \sup \{ \vert \varphi (t_2-t_1,x,\omega)-x\vert_{E} \ : \ \vert x
\vert _{\mathcal D (A)}, \
\omega \in \Omega \} \to 0  \nonumber
\end{eqnarray}
as $t_2-t_1 \to 0$. The theorem is proved.
\end{proof}

\begin{example}\label{ex}
Navier-Stokes equations. We consider the two-dimensional Navier-Stokes system
\begin{eqnarray}\label{eq2.9*}
& & u'+q(t)\sum _{i=1}^{2}u_i\partial _i u =\nu \Delta u -\nabla p + \phi (t)  
\\
& & div \ u =0, \ \ u|_{\partial D} =0 , \nonumber
\end{eqnarray}
where $D $ is an open bounded set with smooth
boundary $ \partial D \in C^{2}$.

The functional setting of the problem is well known
\cite{Lad69},\cite{Tem81} . We denote by $H$ and $V$ the closures
of the linear space $\{ u \in C_0^{\infty}(D)^{2} , \ div \ u=0 \}
$ in $L^{2}(D)^{2}$ and $H_0^{1}(D)^{2}$, respectively. Denote by
$P$ the corresponding orthogonal projection $P: L_2(D)^{2} \to H$.
We further set
$$
A:=-\nu P\Delta , \ B(t)(u,v):=q(t)P(\sum _{i=1}^{2}u_i\partial _i v) .
$$

The Stokes operator $A$ is self-adjoint positive with domain
$\mathcal D (A)$  dense in $H$. The inverse operator is compact.
We define the Hilbert spaces $\mathcal D (A^{\alpha}), \ \alpha
\in (0,1] $ as the domains of the powers of $A$ in the standard
way. Furthermore, $V:=\mathcal D (A^{1/2}),$ and $\vert u \vert
_{\mathcal D (A^{1/2})}=\vert  \nabla u\vert.$

Applying $P$ we write (\ref{eq2.9*}) as the evolution equation of the 
following form
\begin{equation}\label{eq2.9**}
u'+Au + \mathcal B (t)(u,u) = \mathcal F (t), \ \mathcal F (t):= P\phi (t) .
\end{equation}
We suppose that $\mathcal F\in C(\mathbb R ,H) \ (X=H)$ and $\mathcal B \in
C(\mathbb R , L^{2}(H,\mathcal D (A^{-\delta})) \ (F=\mathcal D 
(A^{-\delta})).$
Denote by $Y:=C(\mathbb R , H)\times
C(\mathbb R , L^{2}(H,\mathcal D (A^{-\delta}))$ and $(Y,\mathbb R ,\sigma)$
a dynamical system
of translations (Bebutov's dynamical system, see for example, 
\cite{Sch72},\cite{Sch85} and
\cite{Sel}). Let $\Omega :=H(\mathcal B ,\mathcal F)=
\overline{\{(\mathcal B _\tau , \mathcal F _\tau ) | \
\tau \in \mathbb R \}},$ where $\mathcal B_\tau (t):=\mathcal B (t+\tau)\ ($ 
respectively $
\mathcal F_\tau (t):=\mathcal F (t+\tau)) $ for all $t\in \mathbb R$, by bar 
we denote
a closure in the compact-open topology and $(\Omega ,\mathbb R, \sigma )$ be 
a dynamical
system of translations on $\Omega$.

Along with the equation (\ref{eq2.9**}) we consider its $H$-class
\begin{equation}\label{eq2.9***}
u'+Au + \tilde{\mathcal B} (t)(u,u) = \tilde{\mathcal F} (t),
\end{equation}
where $(\tilde{\mathcal B},\tilde{\mathcal F}) \in H(\mathcal B,\mathcal F)$.
Let $B: \Omega \to L^{2}(H,\mathcal D (A^{-\delta}))$ (respectively $f:\Omega 
\to H)$
be a mapping defined by equality
$$
B(\omega)=B(\tilde{\mathcal B},\tilde{\mathcal F}):=
\tilde{\mathcal B (0)} \ (f(\omega)=f(\tilde{\mathcal B},\tilde{\mathcal 
F}):=
\tilde{\mathcal F (0)} ),
$$
where $\omega = (\tilde{\mathcal B},\tilde{\mathcal F})\in \Omega$, then the 
equation
(\ref{eq2.9**}) and its $H$-class can be written in the form (\ref{eq2.9}).

We now set in the notation above $E=\mathcal D (A^{1/2}), \ X=H, \
F=\mathcal D (A^{-\delta})$ and see that (\ref{eq2.3})-(\ref{eq2.5}), 
(\ref{eq2.6}) and
(\ref{eq5.12*})-(\ref{eq5.13*}) are valid with $\alpha _1=1/2 +\delta,\ \beta 
_1=1/2,\
\beta _2 = 3/2.$
\end{example}

We note that from the conditions (\ref{eq2.7}) -(\ref{eq2.9}) it
follows that
\begin{equation}\label{eq2.10}
\vert B(\omega)(x_1,x_1)-B(\omega)(x_2,x_2)\vert_{F} \le C_B(\vert 
x_1\vert_{E}
+\vert x_2\vert_{E})\vert x_1 -x_2\vert_{E}
\end{equation}
for all $x_1,x_2 \in E$ and $\omega \in \Omega.$

According to Theorem \ref{t*} through every point $ x \in H$ passes
a unique solution $ \varphi (t,x,\omega)$ of equation
(\ref{eq2.1}) at the initial moment $ t=0$. And this solution is
defined on some interval $ [0,t_{(x,\omega)} )$. Let us note, that
\begin{equation}\label{eq2.11}
\begin{array}{ll}
w'(t)\!\!&=2 Re \langle \varphi ' (t,x,\omega ), \varphi (t,x, \omega ) 
\rangle =
      2Re \langle A(\omega t) \varphi (t,x,\omega), \varphi (t,x,\omega) 
\rangle + \\
     & 2Re \langle B(\omega t)(\varphi(t,x,\omega),
     \varphi(t,x,\omega) ), \varphi (t,x,y) \rangle +
       2Re \langle f(\omega t), \varphi (t,x,\omega) \rangle \\
     & =2Re \langle  A(\omega t) \varphi (t,x,\omega), \varphi (t,x,\omega) 
\rangle +
       2Re \langle  f(\omega t), \varphi (t,x,\omega) \rangle \\
      & \le -2 \alpha \vert \varphi (t,x,\omega) \vert ^2_{E} +
 2 \Vert f \Vert \vert \varphi (t,x,\omega) \vert_{E} ,
 \end{array}
\end{equation}
where
$ \Vert f \Vert := \max \{ \vert f(\omega) \vert_{X} : \omega \in \Omega \}$ 
and
$ w(t)= \vert \varphi (t,x,\omega) \vert ^2_{E} .$ Then
\begin{equation}\label{eq2.12}
w' \le -2\alpha w + 2 \Vert f \Vert  w^{\frac{1}{2}}
\end{equation}
and consequently
\begin{equation}\label{eq2.13}
w(t) \le v(t)
\end{equation}
for all $ t \in [0, t_{(x,\omega)} ) $, where $ v(t)$ is an upper solution
of equation
\begin{equation}\label{eq2.14}
v' = -2\alpha v + 2 \vert \vert f \vert \vert v^{\frac{1}{2}},
\end{equation}
satisfying condition $ v(0)=w(0)=\vert x \vert^2 .$ Thus
\begin{equation}\label{eq2.15}
v(t)=[( \vert x \vert_{E} - \frac {\Vert f \Vert }{\alpha} )
e^{-\alpha t} + \frac {\Vert f \Vert } {\alpha} ]^2
\end{equation}
and consequently
\begin{equation}\label{eq2.16*}
\vert \varphi (t,x,\omega) \vert_{E} \le ( \vert x \vert _{E} - \frac{\Vert f
\Vert }{\alpha})e^{-\alpha t} + \frac{\Vert f \Vert }{\alpha}
\end{equation}
for all $ t \in [0, t_{(x,\omega)}).$ It follows from the inequality
(\ref{eq2.12}) that solution $ \varphi (t,x,\omega) $ is bounded and
therefore it may be prolonged on $ \mathbb R_+ = [0,+ \infty ).$

Thus we have proved the following theorem.

\begin{theorem}\label{t2.2}
Let the conditions (\ref{eq2.6}) and (\ref{eq2.7}) are fulfilled. Then
the following statements hold:
\begin{enumerate}
\item[(i)] Every solution $\varphi (t,x,\omega)$ of non-autonomous
Navier-Stokes equation (\ref{eq2.1}) is bounded and therefore it
may be prolonged on $\mathbb R_{+}$. \item[(ii)]
\begin{equation}\label{eq2.17*}
\vert \varphi (t,x,\omega)\vert _{E}\le C(\vert x\vert _{E}),
\end{equation}
for all $\ t\ge 0,\ \omega \in \Omega$ and $x\in E,$
where $C(r)=r$ if $r\ge r_0:=\frac{\Vert f\Vert}{\alpha}$
and $C(r)=r_0$ if $r\le r_0$;
\item[(iii)]
\begin{equation}\label{eq2.13*}
\limsup \limits _{t\to+\infty}\sup\{ \vert\varphi (t,x,\omega) \vert _{E}\ :
\vert x \vert _{E} \le r, \omega \in \Omega\}\le
\frac{\vert\vert f \vert\vert}{\alpha}
\end{equation}
for every $ r>0.$
\end{enumerate}
\end{theorem}

\begin{lemma}\label{l2.10}
Under the conditions of Theorem \ref{t2.2} we have
\begin{equation}\label{eq2.14*}
\int _t^{t+l}\vert  \varphi (\tau ,x,\omega)\vert _{E}^{2}d\tau
\le \frac{r^{2}}{2\alpha}+ \frac{r}{\alpha}l\Vert f \Vert :=M(r)
\end{equation}
for all $t\ge 0$ and $r\ge r_0.$
\end{lemma}
\begin{proof} From the equality (\ref{eq2.12}) after integration in
$t$ between $t$ and $t+l$ we obtain
\begin{equation}\label{eq2.15*}
2\alpha \int _t^{t+l}\vert  \varphi (\tau ,x,\omega)\vert
_{E}^{2}d\tau \le \vert \varphi (t,x,\omega)\vert _{E}^{2}
+2rl\Vert f\Vert
\end{equation}
and,consequently,
\begin{equation}\label{eq2.16}
\int _t^{t+l} \vert \varphi (\tau,x,\omega)\vert _E^{2}d\tau \le
\frac{r^2}{2\alpha} + \frac{r}{\alpha}l\Vert f \Vert :=M(r).
\end{equation}
\end{proof}

\begin{lemma}\label{l2.11}(\cite[Ch.3]{Tem})
(The Uniform Gronwall
Lemma). Let $g,\ h,\ y,$ be three positive locally integrable
functions on $]t_0,\infty [$ such that $y'$ is locally integrable
on $]t_0,\infty [$, and which satisfy
\begin{eqnarray}
& & y' \le gy + h \ \text{for} \ t\ge t_0, \nonumber \\
& & \int_t^{t+l}g(s)ds \le a_1, \ \int _t ^{t+l}h(s)ds\le a_2, \ \int _t
^{t+l}y(s)ds \le a_3 \ \text{for} \ t\ge t_0, \nonumber
\end{eqnarray}
where $l,a_1,a_2,a_3,$ are positive constants. Then
\begin{eqnarray}
y(t+l)\le (\frac {a_3}{l}+a_2)e^{a_1} \ \forall \ t\ge t_0.
\nonumber
\end{eqnarray}
\end{lemma}

\begin{theorem}\label{t2.12}
Under the conditions of Theorem \ref{t2.2} if
\begin{eqnarray}\label{eq2.17}
& & \vert  \langle B(\omega)(u,v),w \rangle \vert \le C \vert u
\vert ^{1/2}\vert Au\vert ^{1/2}\vert v \vert _{1/2} \vert w \vert
\\
& & \forall \ u\in \mathcal D (A),\ v \in D(A^{1/2}),\ w \in X ,
\end{eqnarray}
then
\begin{eqnarray}\label{eq2.18}
\vert \varphi (t,x,\omega) \vert _{\mathcal D (A)} \le K(r) \
\forall \vert x \vert \le r \ (r\ge r_0),
\end{eqnarray}
where $K(r)$ is some positive constant depending only on $r$.
\end{theorem}
\begin{proof} Since
\begin{equation}\label{eq2.19}
\langle A \varphi (t,x,\omega) , \varphi (t,x,\omega) \rangle
=\frac{1}{2}\frac{d}{dt}\vert \varphi (t,x,\omega)\vert _{E}^{2}
\end{equation}
by taking the scalar product of (\ref{eq2.1}) with $Au$ we find
\begin{eqnarray}\label{eq2.20}
& & \frac{1}{2}\frac{d}{dt}\vert \varphi (t,x,\omega)\vert _{E}^{2} +
\vert  A \varphi (t,x,\omega)\vert _{E}^{2} +\nonumber \\
& & \langle B(\varphi (t,x,\omega),
\varphi (t,x,\omega)) , A \varphi (t,x,\omega)\rangle =
\langle f(\omega t), A \varphi (t,x,\omega)\rangle .
\end{eqnarray}
Taking into account the inequality
\begin{equation}\label{eq2.21}
\vert \langle f(\omega t), A \varphi (t,x,\omega)\rangle \vert \le
\vert f(\omega t) \vert _{E} \vert A \varphi (t,x,\omega) \vert _{E} \le
\frac{1}{4}\vert A \varphi (t,x,\omega) \vert ^{2} + \Vert f \Vert ^{2}
\end{equation}
and using (\ref{eq2.19}) and the Young inequality we obtain
\begin{eqnarray}\label{eq2.22}
&&\vert \langle B(\omega)(\varphi (t,x,\omega),\varphi (t,x,\omega)),
A\varphi (t,x,\omega) \rangle \vert \le \nonumber \\
&&c_1\vert \varphi (t,x,\omega) \vert ^{1/2}
\Vert \varphi (t,x,\omega) \Vert
\vert A\varphi (t,x,\omega) \vert^{3/2} \le \\
&&\frac{1}{4}\vert A \varphi (t,x,\omega) \vert ^{2}+
c_1' \vert \varphi (t,x,\omega) \vert ^{2}
\Vert \varphi (t,x,\omega) \Vert ^{4} . \nonumber
\end{eqnarray}
Hence
\begin{eqnarray}\label{eq2.23}
&&\frac{1}{2}\frac{d}{dt}\Vert \varphi (t,x,\omega) \Vert ^{2} +
\vert A \varphi (t,x,\omega) \vert ^{2} \le
\Vert f \Vert ^{2} +\frac{1}{2} \vert A \varphi (t,x,\omega) \vert ^{2} +\\
&&c_1'\vert \varphi (t,x,\omega) \vert ^{2}
\Vert \varphi (t,x,\omega) \Vert ^{4}. \nonumber
\end{eqnarray}

From this inequality according to Gronwal lemma we can prove that
$\vert \varphi (t,x,\omega) \vert _{\mathcal D (A)}$ is uniformly (w.r.t. $x$ 
and
$\omega$) bounded on interval $[0,l].$ Applying the uniform Gronwal lemma 
with
$g, h, y$ replaced by
\begin{equation}\label{eq2.24}
c_1' \vert \varphi (t,x,\omega) \vert ^{2}
\Vert \varphi (t,x,\omega) \Vert ^{2} \ , \ \Vert f \Vert ^{2} \ ,
\ \Vert \varphi (t,x,\omega) \Vert ^{2}
\end{equation}
we obtain that $ \Vert \varphi (t,x,\omega) \Vert ^{2} $ is
bounded on $[l,\infty [$ and, consequently, it is bounded on
$[0,\infty [$ uniformly w.r.t. $\Vert x \Vert \le r$ and $\omega
\in \Omega .$ The theorem is proved.
\end{proof}

\section{Non-autonomous dissipative dynamical systems
and their attractors.}

Let $\Omega $ and $W$ be two metric spaces,
$ (\Omega , \mathbb R ,\sigma )$ be a dynamical system
on $ \Omega .$ Let us consider a continuous mapping
$ \varphi : \mathbb R ^{+} \times \Omega \times W \to W $ satisfying
the following conditions:
$$   \varphi (0,\cdot,\omega )=id_{W}\ \ \varphi (t+\tau ,x,\omega)
=\varphi (t,\varphi (\tau,x,\omega ),\omega \tau) $$
for all $t, \tau \in \mathbb R ^{+}$, $\omega \in \Omega $ and $ x \in W $. 
Such
mapping $ \varphi $ ( or more explicit $\langle W, \varphi ,
(\Omega ,\mathbb  R ,\sigma ) \rangle $) is called  \cite{Arn}, \cite{Sel}
a cocycle on $(\Omega, \mathbb R,\sigma ) $ with fiber $W$.

\begin{remark}\label{r*}
The non-autonomous Navier-Stokes equation (\ref{eq2.1}) generates
a cocycle $ \varphi $ ( or more explicit $\langle E, \varphi ,
(\Omega ,\mathbb  R ,\sigma ) \rangle $), where $\varphi
(t,x,\omega)$ is a unique solution of equation (\ref{eq2.1})
defined on $\mathbb R _{+}$ with the initial condition $\varphi
(0,x,\omega)=x. $

In fact, according to Theorems \ref{t*} the mapping $\varphi :
\mathbb \times E \times \Omega \to E \ ((t,x,\omega) \to \varphi
(t,x,\omega))$ is continuous and in view of uniqueness of solution
$\varphi (t,x,\omega)$ we have the following identity: $\varphi
(t+\tau,x,\omega)=\varphi (t,\varphi (\tau,x,\omega),\omega \tau
)$ for all $t,\tau \in \mathbb R _{+}, x\in E$ and $\omega \in
\Omega$, where $\omega \tau :=\sigma (\tau, \omega)$.
\end{remark}

\begin{example}\label{ex3.1}
Let $E$ be a Banach space and $C(\mathbb R \times E,E)$ be a space of all
continuous functions $F: \mathbb R \times E \to E$ equipped by the 
compact-open
topology. Let us consider a parameterized differential equation
$$ \frac{dx}{dt}+ Ax = F(\sigma _t \omega ,x)   \ \ ( \omega \in \Omega ) $$
on a Banach space $E$ with $\Omega$ $=$ $C(\mathbb{R}\times E,
E),$ where $\sigma _t\omega := \sigma (t,\omega)$ and
the linear operator $A$ is densely defined in
$E$ and such that the linear equation
\begin{eqnarray}\label{eq2.2*}
&& u'+ Au =0   \nonumber
\end{eqnarray}
generates the $c_0$-semigroup of linear bounded operators
$$
e^{-At}: E \to E, \  \varphi (t,x):=e^{-At}x .
$$
We will define
$\sigma_t$ $:$ $\Omega $ $\to$ $\Omega$
by $\sigma_t \omega (\cdot ,\cdot )$ $=$ $\omega (t+\cdot ,\cdot )$
for each $t$ $\in$ $\mathbb{R}$ and
interpret $\varphi(t,x,\omega)$ as mild solution of the initial value problem
\begin{equation} \label{eq3.1}
\frac{d }{dt} x(t) +Ax =  F(\sigma _t \omega ,x(t)), \quad x(0) = x.
\end{equation}
Under appropriate assumptions on $F$ $:$  $\Omega \times E$ $\to$
$E$ (or even
$F$ $:$  $\mathbb{R}$ $\times E$
$\to$ $E $ with $\omega (t)$ instead of $\sigma_t \omega $
in (\ref{eq3.1})) to ensure forwards existence and uniqueness,
then $\varphi $ is a cocycle on
$( C(\mathbb R \times E, E), \mathbb R , \sigma) $ with
fiber $E$, where $( C(\mathbb R \times E, E), \mathbb R , \sigma) $
is a Bebutov's dynamical system
(see for example \cite{Bro79},\cite{Ch02}, \cite{Sch72},\cite{Sel}).
\end{example}

The triplet $\langle (X,\mathbb R _{+},\pi),(\Omega ,\mathbb R
,\sigma),h\rangle $, where $h:X \to \Omega $ is a homomorphism
from the dynamical system $(X,\mathbb R _{+},\pi) $ onto $(\Omega
,\mathbb R ,\sigma)$, is called (see \cite{Bro84},\cite{Ch02}) a
non-autonomous dynamical system.

Let $\varphi $ be a cocycle on $(\Omega ,\mathbb R ,\sigma )$ with the fiber
$E$. Then the mapping $\pi$  $:$
$\mathbb R^{+}$  $\times \Omega
  \times  E$  $\to$  $\Omega \times  E $ defined by
$$
\pi(t,x,\omega) := (\varphi(t,x,\omega), \sigma _t \omega)
$$
for all $t$  $\in$  $\mathbb{R}^{+}$ and  $(x,\omega)$$\in$ $E\times
\Omega $ forms a semi-group
$\{\pi(t,\cdot,\cdot)\}_{t \in\mathbb{R}^{+}}$ of mappings of
$X:=\Omega \times E$
into itself, thus a semi-dynamical system on the
state space
$X$,  which is called a skew-product flow
\cite{Sel} and the triplet
$\langle (X,\mathbb R _{+},\pi),(\Omega ,\mathbb R ,\sigma),h\rangle $\
(where $h:=pr_2 : X \to \Omega )$
 is a non-autonomous dynamical system.

The cocycle $\varphi $ over $(\Omega ,\mathbb R,\sigma)$ with the fiber $W$
we will define by a compact (bounded) dissipative one , if there is a
nonempty  compact $K \subseteq W$ such that
\begin{equation}\label{eq3.2}
 \lim_{t \to + \infty} \sup \{ \beta (U(t,\omega)M,K) \vert
  \omega \in \Omega \}=0
\end{equation}
for any $M \in C(W)$ (respectively $M\in \mathcal B(W)),$ where by
$C(W) (\ \mathcal B (W))$ is denoted a family of all compact (bounded)
 subsets of $W$, by $\beta $ a semi-distance of Hausdorff and $U(t,\omega):=
\varphi (t,\cdot,\omega)$.

\begin{lemma}\label{l3.2} Let $\Omega$ be a
compact metric space and
$\langle W,\varphi,(\Omega ,\mathbb R ,\sigma)\rangle $
be a cocycle over $(\Omega ,\mathbb R ,\sigma)$ with the fiber  $W$.
In order to $\langle W,\varphi,(\Omega ,\mathbb R ,\sigma)\rangle $
be a compact (bounded) dissipative, it is necessary and sufficient
that the skew-product
dynamical system $(X,\mathbb R_+,\pi)$ should be a compact (bounded)
dissipative one .
\end{lemma}

This assertion
directly follows from
the corresponding definitions (see for example \cite{Hal},\cite{Ch02}).

By the whole trajectory of the semi-group dynamical system
$(X,\mathbb R_+,\pi)$ (of the cocycle
$\langle W,\varphi,(\Omega ,\mathbb R ,\sigma) \rangle $ over
$(\Omega ,\mathbb R ,\sigma )$
with the fiber $W$), which passes through the point $x \in X
((u,y) \in W \times \Omega )$ we will call the continuous mapping
$ \gamma : \mathbb R \to X ( \nu : \mathbb R \to W)$ which satisfies
the conditions :  $ \gamma (0)=x \ ( \nu (0)=u) $ and
$ \pi^{t} \gamma (\tau)=\gamma (t+ \tau)
\ (\nu(t+ \tau)=\varphi(t,\nu (\tau),\omega t) )$
for all  $t \in \mathbb R_+$ and  $\tau \in \mathbb R$.
If $ M \subseteq W $, then we denote by
\begin{eqnarray}\label{eq*}
&&\Omega_{\omega}(M)=\bigcap_{t\ge 0} \overline{\bigcup_{\tau\ge t}
\varphi(\tau,M,\omega ^{-\tau})} \nonumber
\end{eqnarray}
for every  $\omega \in \Omega $, where $\omega ^{-\tau}:=
\sigma (-\tau, \omega) $.

\begin{theorem}\label{t3.3} (\cite{Ch97},\cite{Ch02})
Let $\Omega $ be a compact metric space,
$\langle W,\varphi,(\Omega ,\mathbb R ,\sigma)\rangle $
 be a compactly (boundedly) dissipative cocycle and
$K$ be a nonempty compact,
arising in the equality (\ref{eq3.2}), then the following assertions hold:
\begin{enumerate}
\item
$I_\omega := \Omega_{\omega}(K) \ne \emptyset $, is compact,
$ I_{\omega} \subseteq K$
   and
\begin{eqnarray}\label{eq**}
&& \lim_{t \to + \infty} \beta(U(t,\omega ^{-t})K,I_{\omega})=0 \nonumber
\end{eqnarray}
   for every $\omega \in \Omega$;
\item
$U(t,\omega)I_\omega =I_{\omega t} $ for all
$\omega \in \Omega $ and  $t \in \mathbb R_+$;
\item
\begin{eqnarray}\label{eq3.3}
&& \lim_{t \to + \infty}\beta(U(t,\omega ^{-t})M,I_{\omega})=0 \nonumber
\end{eqnarray}
for all  $M \in C(W)$ (respectively $M\in \mathcal B(X)$) and
$\omega \in \Omega $ ;
\item
\begin{eqnarray}\label{eq3.4}
&& \lim _{t \to + \infty} \sup \{ \beta(U(t,\omega ^{-t})M,I) \vert
   \omega \in \Omega \}=0 \nonumber
\end{eqnarray}
for any  $M \in C(W)$ (respectively $M\in \mathcal B(X)$) , where
   $I=\cup \{ I_\omega \vert \ \omega \in \Omega \} $;
\item $ I_\omega :=pr_1 J_\omega $ for all  $\omega \in \Omega $,
where $J$ is a Levinson's center of
   $(X,\mathbb R _+,\pi)$, and, hence, $I=pr_1J$;
\item
the set $I$ is compact;
\item
the set $I$ is connected if the spaces $W$ and $Y$ are connected.
\end{enumerate}
\end{theorem}

The family of compact
sets $\{I_\omega| \omega \in \Omega \} \quad (I_\omega \subset W $
is nonempty compact for every
$\omega \in \Omega )$ is called (see, for example, \cite{Ch97} or 
\cite{Ch02})
the compact global attractor of cocycle $\varphi$ if the following conditions 
are
fulfilled:
\begin{enumerate}
\item
The set $I:=\bigcup \{I_\omega |\quad  \omega \in \Omega \} $ is precompact.
\item
$\{I_\omega |\quad  \omega \in \Omega \}$ is invariant w.r.t. the cocycle
$\varphi ,$ i.e.
$
\varphi (t,\omega,I_\omega)=I_{\sigma_t\omega}
$
for all $t\in \mathbb R _{+}$ and $\omega \in \Omega.$
\item
The equality
$
\lim \limits _{t \to + \infty}\sup \limits _{\omega \in \ \Omega}\beta 
(\varphi (t,K,\omega),I)=0
$
holds for every bounded set $K\subset W $.
\end{enumerate}

\begin{coro}\label{cor3.5} Under the conditions of Theorem \ref{t3.3} the 
cocycle $\varphi$ admits
a compact global attractor.
\end{coro}

Dynamical system $(X, \mathbb R_{+}, \pi )$ is called asymptotically
compact (see \cite{Hal},\cite{Lad} and also
\cite{Ch97},\cite{Ch02}) if for any positive invariant bounded set
$A \subset X$
there is a compact $K_A \subset X$ such that
\begin{equation}\label{eq3.5}
\lim \limits _{t \to + \infty} \beta (\pi^{t}A, K_A)=0 .
\end{equation}

Dynamical system $(X, \mathbb R_{+}, \pi )$ is called compact
(completely continuous)
if for every bounded set $A\subset X$ there exists a positive
number $l=l(A)$ such that the
set $\pi ^{l}A$ will be precompact.

It is easy to verify (see for example \cite{Ch02}) that every
compact dynamical system $(X, \mathbb R_{+}, \pi )$ is
asymptotically compact.

The cocycle $\langle W, \varphi , (Y, \mathbb R, \sigma )\rangle $
is called compact (asymptotically compact)
if the skew-product dynamical system
$(X,\mathbb R_{+},\pi) \ (X=W\times Y, \pi=(\varphi ,\sigma))$ is compact
(respectively asymptotic compact).

Let $(X,\mathbb R _{+},\pi)$ be compact dissipative and  $K$ be a
compact set, which attracts all compact subsets of $X$. Suppose
\begin{equation}\label{eq3.6}
J=\Omega(K),
\end{equation}
where  $\Omega(K)=\bigcap_{t\geq0} \overline {\bigcup_{\tau\geq
t}\pi^\tau K}$. The set $J$ defined by the equality  (\ref{eq3.6})
does not depend on selection of the attractor $K$, and is
characterized only by the properties of the dynamical system
$(X,\mathbb R_+,\pi)$ itself. The set $J$ is called the Levinson's
center of the compact dissipative system $(X,\mathbb R_{+},\pi)$.

\begin{theorem}\label{t3.4} (\cite{Ch97},\cite{Ch02})
Let $(E,\Omega ,h)$ be a local-trivial Banach fibering, $\langle
(E,\mathbb R _{+},\pi),$ $(\Omega ,\mathbb R,\sigma),h \rangle $ be a
non-autonomous dynamical  system and the dynamical system \\
$(E,\mathbb R _{+},\pi)$ be completely continuous. Then next
conditions are equivalent :
\begin{enumerate}
\item
there is a positive number $r$
such that for any $x \in X$
there will be $\tau=\tau (x) \ge 0 $
for which  $\vert x \tau \vert < r $;
\item
the dynamical system $ (E,\mathbb R _{+},\pi) $ is compact
dissipative and
\begin{eqnarray}\label{eq3.7}
&& \lim\limits_{t \to+ \infty} \sup\limits_{\vert x \vert \le R}
{ \rho (xt,J) }=0 \nonumber
\end{eqnarray}
for any  $R > 0$, where  $J$ is a Levinson's center of dynamical
system $(E,\mathbb R _{+},\pi)$, that is the non-autonomous system
$\langle (E,\mathbb R _{+},\pi),(\Omega ,\mathbb R,\sigma),h
\rangle  $ admits a compact global attractor $J$.
\end{enumerate}
\end{theorem}

A dynamical system $(X, \mathbb R _{+},\pi )$ satisfies
the condition of Ladyzhenskaya
(see \cite{Lad} and also \cite{Ch02}) if for any bounded set
$A \subset X$
there is a compact $K_A \subset X$ such that the equality (\ref{eq3.5}) 
holds.

\begin{theorem}\label{t3.5} (\cite{Ch97},\cite{Ch02})
Let $\langle (E,\mathbb R _{+},\pi),(\Omega ,\mathbb R,\sigma),h
\rangle  $ be a non-autonomous
 dynamical system and let $(E,\mathbb R _{+},\pi)$ satisfy the condition of
 Ladyzhenskaya. Then the conditions  1. and 2. of the theorem \ref{t3.4}
 are equivalent.
 \end{theorem}

Applying general theorems about non-autonomous dissipative systems
to non-autonomous system constructed in the example \ref{ex3.1},
we will obtain series of facts concerning the equation
(\ref{eq2.1}). In particular, from Theorems \ref{t2.2}, \ref{t3.3}
and \ref{t3.5} follows the theorem below.

\begin{theorem}\label{t3.6}
Let $ \Omega $ be a compact metric space, $ (\Omega ,\mathbb R ,
\sigma ) $ be a dynamical system on $ \Omega $ and the conditions
(\ref{eq2.6}) and (\ref{eq2.7}) are fulfilled. If the cocycle
$\varphi$ generated by by non-autonomous Navier-Stokes equation is
asymptotically compact, then for every $ \omega \in \Omega $ there
exists a non-empty compact and connected $ I_{\omega} \subset E $
such that the following conditions hold:
\begin{enumerate}
\item the set $ I= \cup \{ I_\omega \ :\  \omega \in \Omega \}$ is
compact in $ E;$ \item $I$ is connected if $\Omega $ is connected;
\item
$$ \lim \limits _{t \to + \infty} \sup \limits _{\omega \in \Omega} \beta
(U(t,\omega ^{-t})M, I)=0 $$ for any bounded set
$ M \subset E $, where
$ U(t,\omega)= \varphi(t, \cdot ,\omega )$
and $ \beta $ is the semi-distance of Hausdorff;
\item
$U(t,\omega)I_\omega =I_{\omega t} $ for all $ t \in \mathbb R_+$ and
$\omega \in \Omega $;
\item
$ I_\omega $ consists of those and only those points $ x \in E $
through which passes the solution of the
equation (\ref{eq2.1}) bounded on $ \mathbb R $.
\end{enumerate}
\end{theorem}

\begin{theorem}\label{t3.7}
Under conditions of Theorem \ref{t3.6}
\begin{eqnarray}\label{eq3.8}
&& \vert \varphi(t,x,\omega) \vert \le
\frac {\Vert f \Vert } {\alpha} \nonumber
\end{eqnarray}
for all $ t \in \mathbb R, \ \omega \in \Omega $ and $ x \in
I_\omega $, where $\varphi$ is a cocycle, generated by
non-autonomous Navier-Stokes equation.
\end{theorem}
\begin{proof}
According to Theorem \ref{t3.3} the set $ J=\bigcup \{ I_\omega
\times \{\omega \} : \omega \in \Omega \} $ is a Levinson's center
of dynamical system $ (X, \mathbb R_+, \pi) $ and according to
(\ref{eq3.6}) for any point $ (u_0,y_0)=z \in J $ there exists $
t_n \to + \infty , u_n \in E $ and $ \omega _n \in \Omega $ such
that  the sequence $\{u_n\}$ is bounded, $ u_0 = \lim \limits_{n
\to + \infty} \varphi(t_n,u_n,\omega_n) $ and $ \omega_0 = \lim
\limits_{n \to + \infty} \omega_n t_n .$ From the inequality
(\ref{eq2.12}) follows that $ \vert u_0 \vert \le \frac {\Vert f
\Vert}{\alpha}, $ i.e. $ \varphi (t,x,\omega) \in I_{\omega t} $
for all $ \omega \in \Omega $ and $ t \ge 0 $, hence $ \vert
\varphi (t,x,\omega) \vert \le \frac {\Vert f \Vert }{\alpha} $
for any $ t \in \mathbb R, x \in I_\omega $ and $ \omega \in
\Omega .$ The theorem is proved.
\end{proof}

\section{Almost periodic and recurrent solutions of non-autonomous
Navier-Stokes equations}

Let $\mathbb T = \mathbb R $ or $\mathbb R_{+}$ and $(X,\mathbb T,\pi)$
be a dynamical system.

\begin{definition}
The point $x\in X$ is called a stationary
($\tau$-periodic, $\tau >0, \tau \in \mathbb T$) point,
if $xt=x$ ($x\tau = x$ respectively)
for all $t\in \mathbb T$, where $xt:=\pi (t,x)$.
\end{definition}

\begin{definition}
The number $\tau \in \mathbb T$ is called $\varepsilon >0$ shift
(almost period) of point $x \in X$ if $\rho (x\tau,x)<\varepsilon
$ (respectively $\rho (x(\tau +t),xt)<\varepsilon$ for all $t\in \mathbb
T$).
\end{definition}

\begin{definition}
The point $x \in X $ is called almost recurrent (almost periodic)
if for any $\varepsilon$ there exists a positive number $l$ such
that on any segment of length $l$, will be found a $\varepsilon$
shift (almost period) of point $x\in X$.
\end{definition}

\begin{definition}
If a point $x\in X$ is almost recurrent and the set
$H(x)=\overline{\{xt\ \vert \ t\in \mathbb T\}}$ is compact,
then $x$ is called recurrent.
\end{definition}

\begin{definition}
An autonomous dynamical system $(\Omega, {\mathbb T},\sigma)$
 is said to be
pseudo recurrent   if the following conditions are fulfilled:
\begin{enumerate}
\item[a)]
$\Omega$  is compact;
\item[b)]
$(\Omega, {\mathbb T},\sigma)$ is transitive, i.e. there exists a point
$\omega_0\in\Omega$ such that
$\Omega=\overline{\{\omega_0t\mid t\in {\mathbb T}\}}$;
\item[c)]
every point $\omega\in\Omega$ is stable in the sense of Poisson, i.e.
\begin{displaymath}
{\frak N}_{\omega}=\{\{t_n\}\mid \omega t_n\to \omega\; {\mbox and}\;
\vert t_n\vert\to +\infty\}\ne\emptyset.
\end{displaymath}
\end{enumerate}
\end{definition}

\begin{definition} A point $x\in X$ is said to be pseudo recurrent is
the dynamical system $(H(x),\mathbb T,\pi)$ is pseudo recurrent.
\end{definition}

\begin{lemma}\label{l6.1} (\cite{CDG})
Let $\langle (X,{\mathbb T}_1,\pi),(\Omega,{\mathbb T}_2,\sigma),h
\rangle $ be a non-autonomous dynamical system and the following
conditions are fulfilled: \begin{enumerate} \item[$1)$]
$(\Omega,{\mathbb T}_2,\sigma)$ is pseudo recurrent; \item[$2)$]
$\gamma\in C(\Omega,X)$ is an invariant section of the
homomorphism $h:X\to\Omega ,$ i.e. $h(\gamma (\omega))=\omega $
and $\gamma (\sigma (t,\omega))=\pi (t,\gamma (\omega))$ for all
$\omega \in \Omega $ and $t\in \mathbb T_2$.
\end{enumerate}

Then the autonomous dynamical system
$(\gamma(\Omega),{\mathbb T}_2,\pi)$ is  pseudo recurrent.
\end{lemma}

Let $\mathbb T=\mathbb S$ and $(X,\mathbb S,\pi)$ be a bi-sided
dynamical system.

\begin{definition}\label{def4.5*}
A recurrent point $x\in X$ is called almost automorphic (see, for
example, \cite{SheYi}) if whenever $t_{\alpha}$ is a net with
$xt_{\alpha}\to x_{*}$, then $x_{*}(-t_{\alpha})\to x$ too.
\end{definition}

\begin{definition}\label{def4.6*}
A motion $\varphi (t,u_0,y_0)$ ($u_0\in E$ and $y_0\in Y$) of the
cocycle $\varphi$ is called recurrent (almost periodic, almost
automorphic, quasi-periodic, periodic), if the point
$x_0:=(u_0,y_0) \in X:=E\times Y$ is a recurrent (almost periodic,
almost automorphic, quasi-periodic, periodic) point of the
skew-product dynamical system $(X,\mathbb S_{+},\pi)$ \ ( $\pi :=
(\varphi, \sigma)$).
\end{definition}

\begin{remark}\label{r4.1*}
We note (see, for example, \cite{Lev-Zhi}, \cite{SacSel} and
\cite{Sch72,Sch85}) that if $y \in Y $ is a stationary
(respectively $\tau$-periodic, almost periodic, quasi periodic, recurrent)
point of the dynamical system $(Y ,\mathbb T_{2} ,\sigma)$ and
$h:Y \to X $ is a homomorphism of the dynamical system $(Y,\mathbb
T _{2},\sigma)$ onto $(X,\mathbb T_{1},\pi)$, then the point
$x=h(y)$ is a stationary (respectively $\tau$-periodic, almost periodic, 
quasi
periodic, recurrent) point of the system $(X,\mathbb T_{1},\pi)$.
\end{remark}

\begin{lemma}\label{l4.1*}
If $y \in Y $ is an almost automorphic point of the dynamical
system $(Y ,\mathbb S ,\sigma)$ and $h:Y \to X $ is a homomorphism
of the dynamical system $(Y,\mathbb S,\sigma)$ onto $(X,\mathbb
S_{+},\pi)$, then the point $x=h(y)$ is an almost automorphic
point of the system $(X,\mathbb S_{+},\pi)$.
\end{lemma}
\begin{proof}
Let $t_{\alpha}$ be a net with $xt_{\alpha}\to x_{*}$, then we
have $yt_{\alpha}\to y_{*}$ ($y:=h(x)$ and $y_{*}:=h(x_{*})$).
Since the point $y$ is almost automorphic, then also
$y_{*}(-t_{\alpha})\to y$ and, consequently,
$x_{*}(-t_{\alpha})=h(y _{*}(-t_{\alpha}))\to h(y)=x$. The lemma
is proved.
\end{proof}

\begin{remark}\label{r4.1**}
Let $X:=E\times Y$ and $\pi := (\varphi , \sigma )$. Then mapping
$h:Y \to X$ is a homomorphism of the dynamical system $(Y ,\mathbb
T_{2} ,\sigma)$ onto $(X,\mathbb T_{1},\pi)$ if and only if
$h(y)=(\gamma (y),y )$ for all $y\in Y$, where $\gamma :Y \to E$
is a continuous mapping with the condition that $\gamma (y
t)=\varphi (t,\gamma (y),y )$ for all $y \in Y $ and $t\in \mathbb
T_{2}.$
\end{remark}

\begin{definition}
The solution $\varphi (t,x,\omega)$ of non-autonomous
Navier-Stokes equation (\ref{eq2.1}) is called recurrent (respectively
pseudo recurrent, almost automorphic, almost periodic, quasi periodic),
if the point $(x,\omega) \in H\times \Omega $ is a recurrent
(respectively pseudo recurrent, almost automorphic, almost periodic, quasi
periodic) point of skew-product dynamical system $(X,\mathbb
R_{+},\pi)$ \ ($X=H\times \Omega$ and $\pi = (\varphi, \sigma)$).
\end{definition}

Let $X=H\times \Omega $ and $\pi = (\varphi , \sigma )$, then
mapping $h:\Omega \to X$ is a homomorphism of dynamical system
$(\Omega ,\mathbb R ,\sigma)$ onto $(X,\mathbb R_{+},\pi)$ if and
only if $h(\omega)=(u(\omega),\omega )$ for all $\omega \in \Omega
$, where $u:\Omega \to H$ is a continuous mapping with the
condition that $u(\omega t)=\varphi (t,u(\omega),\omega )$ for all
$\omega \in \Omega $ and $t\in \mathbb R _{+}.$

The following affirmations hold:

\begin{lemma}\label{l4.1}
Let $\Omega $ be a compact metric space, $A$ be a linear operator
densely defined in $E$ such that the equation
\begin{eqnarray}\label{eq4.1*}
&& x'+Ax=0 \nonumber
\end{eqnarray}
generates a $c_0$-semigroup
$\{U(t)\}_{t\ge 0}$. If the condition (\ref{eq2.6}) is fulfilled, then
\begin{eqnarray}\label{eq4.1}
&& \Vert U(t)\Vert \le \exp (-\alpha t) \nonumber
\end{eqnarray}
for all $t\in \mathbb R _{+}$, where
$U(t)$ is a Cauchy's operator of equation (\ref{eq4.1*}).
\end{lemma}
\begin{proof} Let $\varphi (t,x):=U(t)x$,
then according to the inequality \ref{eq2.6} we have
$$
\frac{d}{dt} \vert \varphi (t,x) \vert ^{2} \le -2\alpha \vert
\varphi (t,x) \vert ^{2}
$$
and, consequently,$ \vert
\varphi (t,x) \vert \le \exp (-\alpha t) \vert x \vert$ for
all $x\in H $ and $t\in \mathbb R _{+}.$ Thus
we have $\vert U(t)x\vert \le \exp (-\alpha t) \vert x
\vert $, therefore  $\Vert U(t) \Vert \le \exp (-\alpha t)$
for all $t\in \mathbb R _{+}.$
\end{proof}

\begin{lemma}\label{l4.2}
Suppose that the condition (\ref{eq2.6})
is fulfilled. Then for
every function $f\in C(\Omega ,H)$ there exists a unique function
$\gamma \in C(\Omega ,H)$ defined by equality
\begin{eqnarray}
&& \gamma (\omega )= \int _{-\infty}^
{0}U(-\tau)f(\omega \tau)d\tau  \nonumber
\end{eqnarray}
such that
\begin{eqnarray}\label{eq4.4}
&& \gamma (\omega t)=\varphi (t,\gamma (\omega),\omega)
\end{eqnarray}
for every $\omega \in \Omega $ and $t\in \mathbb R_{+}$,
where $\varphi (t,x,\omega) $ is a solution of equation
\begin{eqnarray}\label{eq4.5}
&& u'=A u+f(\omega t) \nonumber
\end{eqnarray}
with the initial condition $ \varphi (0,x,\omega )=x$ and the following
inequality
\begin{eqnarray}\label{eq4.6}
&& \Vert \gamma \Vert \le \frac{1}{\alpha} \Vert f \Vert  \nonumber
\end{eqnarray}
takes place.
\end{lemma}
\begin{proof} The formulated statement results from Lemma \ref{l4.1} and
Proposition 7.33 from \cite{CL}.
\end{proof}

\begin{lemma}\label{l4.3} Let $\Omega$ be a compact metric space,
the cocycle $\varphi ,$ generated by the non-autonomous
Navier-Stokes equation (\ref{eq2.1}) and $ \alpha ^{-2}\Vert f
\Vert C_B <1 $, then the following inequality
\begin{eqnarray}\label{eq4.7}
& & \vert \varphi (t,x_1,\omega )-\varphi (t,x_2,\omega )\vert \le
e^{-(\alpha -C_B\frac{\Vert f \Vert}{\alpha})t} \vert x_1 - x_2\vert  
\nonumber \\
& & \ (x_1,x_2\in B[0,r_0],
\omega \in \Omega \ \text{and} \ t\in \mathbb R_{+}) \nonumber
\end{eqnarray}
takes place.
\end{lemma}
\begin{proof}
Let $r_0:=\frac{\Vert f \Vert}{\alpha}$ and $x_1,x_2\in B[0,r_0] :=
\{ x \in E \ : \vert x \vert \le r_0 \}.$ According to
Theorem \ref{t2.2} we have $\vert \varphi (t,x_i,\omega )\vert \le r_0 $ for 
all
$t\ge 0, \omega \in \Omega $ and $i=1,2.$ Denote by $ \psi (t) :=
\varphi (t,x_1,\omega )-\varphi (t,x_2,\omega ), $ then we obtain
\begin{eqnarray}\label{eq4.8}
& & \frac{d}{dt}\vert \psi (t)\vert ^{2}=2 Re \langle A\psi (t),\psi 
(t)\rangle +
2 Re \langle B(\omega t)(\psi (t),\varphi (t,x_2,\omega ) ),\psi (t)\rangle 
\le \nonumber \\
& & -2\alpha \vert \psi (t) \vert ^{2} +2 C_B \vert \varphi (t,x_2,\omega )
\vert \vert \psi (t) \vert ^{2} \le \nonumber \\
& & -2\alpha \vert \psi (t) \vert ^{2} +
2C_B \frac{\Vert f \Vert}{\alpha} \vert \psi (t) \vert ^{2} =
-2(\alpha - C_B \frac{\Vert f \Vert}{\alpha})\vert \psi (t) \vert ^{2} 
\nonumber
\end{eqnarray}
and, consequently,
\begin{eqnarray}\label{eq4.8*}
&& \vert \psi (t) \vert ^{2} \le e^{-2(\alpha -
C_B \frac{\Vert f \Vert}{\alpha})t}\vert \psi (0)\vert ^{2} . \nonumber
\end{eqnarray}
Thus we have
\begin{eqnarray}\label{eq4.}
& & \vert \varphi (t,x_1,\omega )-\varphi (t,x_2,\omega )\vert \le
e^{-(\alpha -C_B\frac{\Vert f \Vert}{\alpha})t} \vert x_1 - x_2\vert 
\nonumber \\
& & (x_1,x_2\in B[0,r_0],
\omega \in \Omega \ \text{and} \ t\in \mathbb R_{+}) \nonumber
\end{eqnarray}
for all $ x_1,x_2\in B[0,r_0], \omega \in \Omega $ and $ t\in \mathbb R_{+}$.
The Lemma is proved.
\end{proof}

\begin{theorem}\label{t4.4} Let $r_0:=\frac{\Vert f \Vert}{\alpha}$ and
$\Omega $ be a compact metric space,
the cocycle $\varphi ,$ generated by the non-autonomous
Navier-Stokes equation (\ref{eq2.1}) and $\frac{\Vert f \Vert
C_{B}}{{\alpha}^2}< 1 $. Then there exists a function $\gamma \in
C(\Omega ,B[0,r_0])$ such that:
\begin{enumerate}
\item[a.]
\begin{equation}\label{eq4.9}
\gamma (\omega t)=\varphi (t,\gamma (\omega),\omega)
\end{equation}
for every $\omega \in \Omega$ and $t\in \mathbb R_{+}$,
where $\varphi (t,x,\omega)
$ is a solution of equation (\ref{eq2.1})
with the initial condition $ \varphi (0,x,\omega )=x$;
\item[b.]
\begin{eqnarray}\label{eq4.10}
&& \Vert \gamma \Vert \le \frac{\Vert f \Vert}{\alpha} ; \nonumber \\
\end{eqnarray}
\item[c.]
\begin{eqnarray}\label{eq4.11}
& & \vert \varphi (t,x,\omega )-\gamma (\omega t) \vert \le
e^{-(\alpha -C_B\frac{\Vert f \Vert}{\alpha})t}
\vert x - \gamma (\omega)\vert
\end{eqnarray}
for all $x\in E$, $\omega \in \Omega$ and $t\in \mathbb R_{+}$, where
$\Vert \gamma \Vert :=\sup \{\vert \gamma (\omega)\vert \ :
\omega \in \Omega \}.$
\end{enumerate}
\end{theorem}
\begin{proof}
Let $\Gamma := C(\Omega , B[0,r_0])$ ($C(\Omega ,E)$) be a space
all the continuous functions $f:\Omega \to B[0,r_0]$ (respectively
$f:\Omega \to E$ ) equipped with the distance
\begin{eqnarray}\label{eq4.12}
d(f_1,f_2)=\max \{\vert f_1(\omega)-f_2(\omega)\vert \ : \
\omega \in \Omega \} \ . \nonumber
\end{eqnarray}
Then $(\Gamma ,d) $ ( respectively $(C(\Omega ,E),d) $) is a
complete metric space.

Let $t\in \mathbb R_{+}$. We define the mapping
$ S^{t} : \Gamma \to C(\Omega ,E) $ by the equality
$$
(S^{t}\nu)(\omega):=U(t,\omega ^{-t})\nu (\omega ^{-t})
$$
for all $\omega \in \Omega ,$ where $\omega ^{-t}:=\sigma (-t,\omega)$ and
$U(t,\omega) :=\varphi (t,\cdot ,\omega ).$ According to Theorem \ref{t2.2}
we have $S^{t}(\Gamma)\subseteq \Gamma $ for all $t \in \mathbb R_{+}.$
It easy to see that the family of mappings $\{ S^{t}\ \vert \
t\in \mathbb R_{+}\}$ possesses the following properties:
\begin{enumerate}
\item
$$ S^{0}=Id_{\Gamma} $$
and
\item
$$ S^{t+\tau}=S^{t}S^{\tau} $$
\end{enumerate}
for all $t,\tau \in \mathbb R_{+}.$

Thus $\{ S^{t}\ \vert \ t\in \mathbb R_{+}\}$ forms a commutative
semigroup with identity element. Now we will show that the mapping
$S^{t} ( t>0)$ is a contraction. In fact, let $\nu _1, \nu _2 \in
\Gamma ,$ then we have
\begin{equation}\label{eq4.13}
(S^{t}\nu _1)(\omega)-(S^{t}\nu _2)(\omega)=U(t,\omega ^{-t})
\nu _1(\omega ^{-t})- U(t,\omega ^{-t})\nu _2(\omega ^{-t}).
\end{equation}
From the lemma \ref{l4.3} and the equality (\ref{eq4.13}) it follows that
\begin{eqnarray}\label{eq4.14}
d(S^{t}\nu _1,S^{t}\nu _2) \le
e^{-(\alpha -C_B\frac{\Vert f \Vert}{\alpha})t} d(\nu_1,\nu_2) \nonumber
\end{eqnarray}
for all $t\in \mathbb R_{+}$ and, consequently, there exists a unique common
fixe point $\gamma \in \Gamma ,$ i.e. $ S^{t}\gamma =\gamma$ for all
$t\in \mathbb R_{+}.$ In particularly
$$
U(t,\omega ^{-t})\gamma (\omega ^{-t}) = \gamma (\omega )
$$
for all $t\in \mathbb R_{+}$ and $\omega \in \Omega $. From this equality
follows that
$$
\gamma (\omega t)= U(t,\omega )\gamma (\omega)=
\varphi (t,\gamma (\omega),\omega)
$$
for all $t\in \mathbb R_{+} $ and $\omega \in \Omega .$

Let $ x\in E, \varphi (t,x,\omega) $ be a unique solution of equation
(\ref{eq2.1}) with the initial condition $ \varphi (0,x,\omega)=x $ and
$\gamma \in \Gamma $ the function with the property (\ref{eq4.9}).
Denote by $ \psi (t) :=
\varphi (t,x,\omega )-\gamma (\omega t), $ then we have
\begin{eqnarray}\label{eq4.15}
& & \frac{d}{dt}\vert \psi (t)\vert ^{2}=2 Re \langle A\psi (t),\psi 
(t)\rangle +
2 Re \langle B(\omega t)(\psi (t),\gamma (\omega t) ),\psi (t)\rangle \le 
\nonumber \\
& & -2\alpha \vert \psi (t) \vert ^{2} +
2C_B \vert \gamma (\omega t) \vert \vert \psi (t) \vert ^{2} \le
-2\alpha \vert \psi (t) \vert ^{2} + \nonumber \\
& & 2C_B \frac{\Vert f \Vert}{\alpha} \vert \psi (t) \vert ^{2} =
-2(\alpha - C_B \frac{\Vert f \Vert}{\alpha})\vert \psi (t) \vert ^{2} 
\nonumber
\end{eqnarray}
and, consequently,
\begin{eqnarray}\label{eq4.16}
\vert \psi (t) \vert ^{2} \le e^{-2(\alpha -
C_B \frac{\Vert f \Vert}{\alpha})t}\vert \psi (0)\vert ^{2} . \nonumber
\end{eqnarray}
Thus we have
$$
\vert \varphi (t,x,\omega )-\gamma (\omega t)\vert \le
e^{-(\alpha -C_B\frac{\Vert f \Vert}{\alpha})t} \vert x -
\gamma (\omega)\vert
$$
for all $ x \in E$, $ \omega \in \Omega $ and $ t\in \mathbb R_{+}$.
The theorem is proved.
\end{proof}
\begin{coro}\label{cor4.5}
Under the conditions of Theorem \ref{t4.4} there exists a unique function
$\gamma \in C(\Omega ,E)$ such that
\begin{equation}\label{eq4.17}
\gamma (\omega t)=\varphi (t,
\gamma (\omega),\omega)
\end{equation}
for all $t \in \mathbb R_{+}$ and $\omega \in \Omega .$
\end{coro}
\begin{proof}
Let $\overline{\gamma} \in C(\Omega ,E)$ be a function satisfying the
equality (\ref{eq4.17}) and $\gamma \in \Gamma =
C(\Omega ,B[0,r_0]) $ the function from Theorem \ref{t4.4}. Since
$ \overline{\gamma} (\omega t)=\varphi (t,\overline{\gamma} (\omega),\omega)$
for all $t\in \mathbb R_{+}$ and $\omega \in \Omega ,$ then according to the
inequality (\ref{eq4.11}) we have
\begin{eqnarray}\label{eq4.18}
\vert \overline{\gamma} (\omega t) -\gamma (\omega t) \vert \le
e^{-(\alpha -C_B\frac{\Vert f \Vert}{\alpha})t}
\vert \overline{\gamma}(\omega) - \gamma (\omega)\vert
\end{eqnarray}
for all $t \in \mathbb R_{+}$ and $\omega \in \Omega .$ In particularly, from
(\ref{eq4.18}) we obtain
\begin{eqnarray}\label{eq4.19}
& & \vert \overline{\gamma} (\omega)  -\gamma (\omega) \vert \le
e^{-(\alpha -C_B\frac{\Vert f \Vert}{\alpha})t}
\vert \overline{\gamma}(\omega^{-t}) - \gamma (\omega^{-t})\vert \le \\
& & e^{-(\alpha -C_B\frac{\Vert f \Vert}{\alpha})t}
\Vert \overline{\gamma} - \gamma \Vert \nonumber
\end{eqnarray}
for all $t \in \mathbb R_{+}$ and $\omega \in \Omega ,$ where
$\omega ^{-t}:=\sigma (-t,\omega)$ and
$\Vert \overline{\gamma} - \gamma \Vert :=
\max \{ \vert \overline{\gamma}(\omega)-\gamma (\omega)\vert \
: \omega \in \Omega \}.$ Passing to the limit in the inequality 
(\ref{eq4.19})
we obtain $\overline{\gamma}(\omega)=\gamma (\omega) $ for all $\omega \in 
\Omega .$
\end{proof}

\begin{coro}\label{cor4.6}
Under the conditions of Theorem \ref{t4.4}
the equation (\ref{eq2.1}) admits a compact global attractor
$\{I_{\omega} \ : \ \omega \in \Omega \}$ and $I_\omega
= \{\gamma (\omega)\}$ for all $\omega \in \Omega ,$ where
$\gamma \in \Gamma $ is a function from Theorem \ref{t4.4}.
\end{coro}

\begin{coro}\label{cor4.7} The following statements hold:
\begin{enumerate}
\item
Let $\Omega $ be a compact minimal set containing only the
periodic (respectively quasi periodic, almost periodic,
almost automorphic,recurrent) motions, then under conditions of
Theorem \ref{t4.4} the non-autonomous Navier-Stokes equation
(\ref{eq2.1})) admits a unique periodic (respectively
quasi periodic, almost
periodic, almost automorphic, recurrent)
solution $\gamma (\omega t)$ and every other solution of this
equation is asymptotically periodic (respectively asymptotically quasi
periodic, asymptotically almost periodic, asymptotically
automorphic, asymptotically recurrent).
\item
If $(\Omega,\mathbb T,\sigma) $ is a pseudo recurrent dynamical system,
then under conditions of
Theorem \ref{t4.4} the non-autonomous Navier-Stokes equation
(\ref{eq2.1}) admits a unique pseudo recurrent
solution $\gamma (\omega t)$ and every other solution of this
equation is asymptotically pseudo recurrent.
\end{enumerate}
\end{coro}
\begin{proof}
Let $\gamma \in \Gamma $ be a function from Theorem \ref{t4.4},
then according this theorem we have $\varphi (t,\gamma (\omega ),
\omega )=\gamma (\omega t)$ for all $t\in \mathbb R_{+}$ and
$\omega \in \Omega $ and, consequently, the solution $ \varphi
(t,\gamma (\omega ),\omega )$ is periodic (quasi periodic, almost
periodic, almost automorphic, recurrent, pseudo recurrent). Let
$\varphi (t,x,\omega ) $ be a arbitrary solution of equation
(\ref{eq2.1}), then taking into consideration the inequality
(\ref{eq4.11}) we conclude that $\varphi (t,x,\omega ) $ is
asymptotically periodic (asymptotically quasi periodic,
asymptotically almost periodic,asymptotically almost automorphic,
asymptotically recurrent, asymptotically pseudo recurrent).
\end{proof}

\section{Uniform averaging for a finite interval}

We shall be dealing with the non-autonomous Navier-Stokes equation
\begin{equation}\label{eq5.1}
u'+\varepsilon Au+\varepsilon B(u,u)=\varepsilon f(\omega t),
\end{equation}
where $\varepsilon \in [0, \varepsilon _0] $, $A$ is linear and
$B$ is a bilinear operator, $f$ is a forcing term.

Below we will use some notions, denotations and results from \cite{Il98}.
Let Banach spaces $E, F, X, \mathcal E $ satisfy
$$
E\subset F; \ \ \ E, F, X \subset \mathcal E ,
$$
each embedding being dense and continuous.

We suppose that the linear equation
\begin{equation}\label{eq5.2}
u'=Au
\end{equation}
generates the $c_0$-semigroup of linear bounded operators
\begin{equation}\label{eq5.3}
e^{At} : \mathcal E \to \mathcal E ,
\end{equation}
which for $t>0$ can be extended to the linear bounded operators
from $F$ to $E$ satisfying the estimates
(\ref{eq2.3})-(\ref{eq2.5}).

We also suppose that the following condition is satisfied
\begin{equation}\label{eq5.7}
 Ae^{At}=e^{At}A ,
\end{equation}
in the sense of $ L (F,E):=\{A:F \to E \ | A \ \text{is linear and
bounded} \ \}$ equipped with the operational norm.

Function $f$. The external force $f: \Omega \to X $ is continuous, i.e.
$f\in C(\Omega ,X)$.

Operators $e^{At} .$ The operators $e^{-At} \ (t>0)$
can be extended to the linear
bounded operators from $X$ to $E$ satisfying the estimates (\ref{eq5.12*})-
(\ref{eq5.13*}) and the equation (\ref{eq5.7}),
this time in the sense of $\mathcal (X,E)$ .

Existence of partial averaged. $f(\omega )=
f_0(\omega)+f_1(\omega) \ (f_0,f_1 \in C(\Omega,X))$ for all
$\omega \in \Omega $ and the average of $f_1(\omega)$ is equal to
$0$, that is,
\begin{equation}\label{eq5.13}
\lim \limits _{t\to \infty} \frac{1}{t}\int _{0}^{t} f_{1}(\omega \tau )d\tau 
=0
\end{equation}
uniformly with respect to $\omega \in \Omega$.

\begin{remark}\label{r5.2}
1. The condition $(\ref{eq5.13})$ is fulfilled if a dynamical system
 $(\Omega,\mathbb R,\sigma)$ is strictly ergodic, i.e. on $\Omega$
 exists a unique invariant measure $\mu$ w.r.t. $(\Omega,\mathbb R,\sigma)$.

2. According to Lemma 5.1 from \cite{CD} the equality (\ref{eq5.13}) takes 
place
if and only if there exists a positive continuous on $\mathbb R_{+}$
function $k$
with $\lim \limits _{t\to \infty }k(t) =0$ such that
\begin{equation}\label{eq5.14}
\vert \frac{1}{t}\int _{0}^{t} f_{1}(\omega \tau )d\tau \vert_{X} \le k(t)
\end{equation}
for all $\omega \in \Omega $ and $t \in \mathbb R_{+} .$
\end{remark}

Along with equation (\ref{eq5.1}) we consider also the partial
averaged equation
\begin{equation}\label{eq5.16}
u'+\varepsilon
Ax+\varepsilon B(u,u) =
\varepsilon f_0(\omega t).
\end{equation}
If we introduce the "slow time" $\tau:=\varepsilon t$
$(\varepsilon>0)$, then the equations (\ref{eq5.1}) and
(\ref{eq5.16}) can be written in the following way
\begin{equation}\label{eq5.17}
u'+Au +B(u,u)= f(\omega \frac{\tau}{\varepsilon})
\end{equation}
and
\begin{equation}\label{eq5.18}
\bar{u}' + A\bar{u}+ B(\bar{u},\bar{u})=
f_0(\omega \frac{\tau}{\varepsilon}) .
\end{equation}

We will consider the mild solutions $u(t)$ and $\bar{u}(t)$ of the equations
(\ref{eq5.17}) and (\ref{eq5.18}), i.e. $u,\bar{u} \in C([0,T],E)$ and 
satisfy the
following integral equations
\begin{equation}\label{eq5.19}
u(\tau)=e^{-A\tau}x +\int _{0}^{\tau}e^{-A(\tau-s)}
(-B(u(s),u(s))+f(\omega \frac{s}{\varepsilon}))ds ,
\end{equation}
and
\begin{equation}\label{eq5.20}
\bar{u}(\tau)=e^{-A\tau}x +\int _{0}^{\tau}e^{-A(\tau-s)}
(-B(\bar{u}(s),\bar{u}(s))+f_0(\omega \frac{s}{\varepsilon}))ds .
\end{equation}

Denote by $\varphi (\tau ,x,\omega ,\varepsilon) \ (\bar{\varphi}
(\tau ,x,\omega ,\varepsilon) ) $ a unique solution of equation
(\ref{eq5.19}) (respectively (\ref{eq5.20})). According to Theorem
\ref{t2.2} the cocycle $\varphi (\cdot,\cdot,\cdot,\varepsilon) \
(\bar{\varphi}(\cdot,\cdot,\cdot,\varepsilon)) $, generated by
equation (\ref{eq5.19}) (respectively (\ref{eq5.20})), has an
absorbing ball $B[0,R_0] \ (B[0,\bar{R}_0])$ in $E$, where
$R_0:=\frac{\Vert f \Vert}{\alpha} \ (\bar{R}_0:=\frac{\Vert f_0
\Vert}{\alpha})$. This means that for every ball $B[0,R]$
(respectively $B[0,\bar{R}]$) there exists a positive number
$L=L(R) $ (respectively $\bar{L}=\bar{L}(\bar{R})$) such that
\begin{equation}\label{eq5.21}
U(t,\omega,\varepsilon)B[0,R]\subseteq B[0,R_0]
\end{equation}
\begin{equation}\label{eq5.22}
(\bar{U}(t,\omega,\varepsilon)B[0,\bar{R}]\subseteq B[0,\bar{R}_0] )
\end{equation}
for all $t\ge L \ ( t \ge \bar{L}), \varepsilon \in [0,\varepsilon_0]$ and
$\omega \in \Omega $, where $U(t,\omega,\varepsilon)
:=\varphi (t,\cdot ,\omega ,\varepsilon ) $ \ (
$\bar{U}(t,\omega,\varepsilon):
=\bar{\varphi} (t,\cdot ,\omega ,\varepsilon ) $.

According to Theorem \ref{t2.2} the cocycle
$\varphi (\cdot,\cdot,\cdot,\varepsilon) \
(\bar{\varphi}(\cdot,\cdot,\cdot,\varepsilon)) $
is uniformly bounded for $t\ge 0$, that is,
for every ball $B[0,R_1]$\ ($B[0,\bar{R}_1] $) there exists a ball
$B[0,R_2]$\ ($B[0,\bar{R}_2] $) such that
\begin{equation}\label{eq5.23}
U(t,\omega,\varepsilon)B[0,R_1]\subseteq B[0,R_2]
\end{equation}
\begin{equation}\label{eq5.24}
(\bar{U}(t,\omega,\varepsilon)B[0,\bar{R}_1]\subseteq B[0,\bar{R}_2] )
\end{equation}
for all $t\ge 0, \varepsilon \in [0,\varepsilon_0]$ and
$\omega \in \Omega $.

\begin{theorem}\label{t5.3}
Let $L>0$ be arbitrary but fixed. If
$\varphi (0,x,\omega ,\varepsilon )=\bar{\varphi} (0,x,\omega ,\varepsilon )=
x \in B[0,\overline{R}_0] $, that is, the initial points coincide and
belong to the absorbing ball of equation (\ref{eq5.18}) and
the condition (\ref{eq2.17}) is fulfilled, then the following relation
takes place
\begin{equation}\label{eq5.25}
\sup \{ \vert \varphi (t,x,\omega ,\varepsilon )-
\bar{\varphi}(t,x,\omega ,\varepsilon )\vert _{E} \ : \
0\le t \le L, \ \vert x \vert _{E} \le \overline{R}_0 , \ \omega \in
\Omega \} \to 0
\end{equation}
as $\varepsilon \to 0 $.
\end{theorem}
\begin{proof}
The proof below goes along the same lines as the proofs of the
corresponding results from \cite{DF},\cite{Il96} and \cite{Il98}.
We set $v(t):= \varphi (t,x,\omega ,\varepsilon )-
\bar{\varphi}(t,x,\omega ,\varepsilon ) $. Substracting the
equation (\ref{eq5.19}) from the equation (\ref{eq5.20}), we
obtain
\begin{eqnarray}\label{eq5.26}
& & v(t)=\int _0^{t} e^{(t-s)A}(-B(v(s),
\varphi (s,x,\omega,\varepsilon)) - B(\bar{\varphi} 
(s,x,\omega,\varepsilon),v(s)))ds + \nonumber \\
& & \int _0^{t}e^{(t-s)A}f_1(\omega s) ds
\end{eqnarray}

According to Theorem \ref{t2.2} $\vert \varphi (t,x,\omega ,
\varepsilon )\vert , \vert \bar{\varphi}(t,x,\omega ,\varepsilon )
\vert \le r_0 $ for all $t\ge 0,$ where $r_0:=\max \{\frac{\Vert f 
\Vert}{\alpha},
\frac{\Vert f_1 \Vert}{\alpha}\}.$ In view of (\ref{eq5.26}) $v(t)$ satisfies 
the
inequality
\begin{eqnarray}\label{eq5.27}
&& \vert v(t)\vert _{E}\le \vert
\int _0^{t} e^{(t-s)A}(B(v(s),
\varphi (s,x,\omega,\varepsilon)) +
B(\bar{\varphi} (s,x,\omega,\varepsilon),v(s)))ds\vert _{E} + \nonumber \\
& & \vert \int _0^{t}e^{(t-s)A} f_1(\omega \frac{s}{\varepsilon})
ds\vert _{E} \ \nonumber , \ t\in [0,L] .
\end{eqnarray}

By (\ref{eq2.10}) and (\ref{eq2.11}) we see that the first term on the 
right-hand
side of (\ref{eq5.27}) is less than
\begin{equation}\label{eq5.28}
2r_0\cdot K\cdot C_B \int _0^{t} e^{-a(t-s)}(t-s)^{\alpha _1}
\vert v (s)\vert _{E} ds \ .
\end{equation}

We now show that the second term
in (\ref{eq5.27}) tends to
$0$ as $\varepsilon \to 0$ uniformly w.r.t.
$t\in [0,L], \vert x \vert \le R_0 $ and $\omega \in \Omega .$
Integrating by part in $s$ and
taking into account the inequalities (\ref{eq2.4})-(\ref{eq2.5}),
(\ref{eq5.12*})-(\ref{eq5.13*}) and (\ref{eq5.14}) we find
\begin{eqnarray}\label{eq5.28*}
& & \vert \int _0^{t}e^{(t-s)A} f_1(\omega \frac{s}{\varepsilon})
ds\vert _{E} = \vert e^{At}\int _0^{t}f_1(\omega
\frac{s}{\varepsilon}) ds + \int _0^{t} Ae^{(t-s)A}
\int_{t}^{s}f_1(\omega \frac{\tau}{\varepsilon}) d\tau\vert _{E}ds
 \le \nonumber
\\ & & \Vert e^{A\tau} \Vert _{X\to E}\vert \int _0^{t}f_1(\omega
\frac{s}{\varepsilon}) ds \vert _{X} + \int _0^{t}\Vert A
e^{A(t-s)} \Vert _{X\to E} \vert \int _t^{s}f_1(\omega
\frac{s}{\varepsilon}) ds  \vert _{X} \le
\\ & & Kt^{1-\beta _1}e^{-at}k_1(\frac{t}{\varepsilon}) + \int _{0}^{t} K
(t-s)^{1-\beta _2} e^{-a(t-s)} k_1(\frac{t-s}{\varepsilon})ds
.\nonumber
\end{eqnarray}

Let $\alpha \in [0,1), \nu \in (0,1)$ and $\beta \in [0,2)$. Since
\begin{eqnarray}\label{eq5.28**}
&& t^{\alpha}k_1(\frac{t}{\varepsilon})\le \sup \limits_{0\le t
\le \varepsilon ^{\nu}}
t^{\alpha}k_1(\frac{t}{\varepsilon})+\nonumber \sup
\limits_{\varepsilon ^{\nu}\le t \le L} t^{\alpha}k_1(t)\le
\\ && \varepsilon^{\alpha \nu}k_1(0) +L^{\alpha}k_{1}(\varepsilon
^{\nu -1}) :=c_1(\varepsilon) \to 0 \ \text{as} \ \varepsilon \to 0 \nonumber
\end{eqnarray}
and
\begin{eqnarray}\label{eq5.28***}
&&\int_{0}^{t}s^{1-\beta}k_{1}(\frac{s}{\varepsilon})ds
=\int_{0}^{\varepsilon ^{\nu}
}s^{1-\beta}k_{1}(\frac{s}{\varepsilon})ds +\int_{\varepsilon
^{\nu}}^{t}s^{1-\beta}k_{1}(\frac{s}{\varepsilon})ds \le \nonumber
\\ && k_{1}(0)\frac{\varepsilon^{\nu (2-\beta)}}{2-\beta} +
k(\varepsilon ^{\nu -1})\frac{(t^{2-\beta}-\varepsilon^{\nu
(2-\beta)})}{2-\beta} \le \nonumber \\ &&
k_{1}(0)\frac{\varepsilon^{\nu (2-\beta)}}{2-\beta} +
k(\varepsilon ^{\nu -1})\frac{(L^{2-\beta}-\varepsilon^{\nu
(2-\beta)})}{2-\beta} :=c_2(\varepsilon)\ \text{as} \ \varepsilon \to 0
\end{eqnarray}
uniformly w.r.t. $t\in [0,L]$ and $\omega \in \Omega ,$ then
\begin{eqnarray}\label{eq}
\sup \limits _{0\le t \le L}\vert \int _0^{t}e^{(t-s)A} f_1(\omega
\frac{s}{\varepsilon}) ds\vert _{E} \to 0 \ \text{as} \
\varepsilon \to 0.
\end{eqnarray}

Thus,
\begin{equation}\label{eq.}
\vert v(t)\vert _{E} \le C(\varepsilon) + D\cdot \int _0^{t}(t-s)^{-\alpha 
_1}\vert v(s)\vert _{E}ds
\end{equation}
for all $t\in [0,L]$, where $C(\varepsilon):= c_1(\varepsilon)+
+c_2(\varepsilon) \to 0$ as
$\varepsilon \to 0$ and $D:= 2R_0C_BK$.

We now use the known inequality \cite[Ch.7]{Hen}. If
\begin{equation}\label{eq..}
u(t) \le a + b\int _0^{t}(t-s)^{\beta -1}u(s)ds, \ 0<\beta \le 1 ,
\end{equation}
then
\begin{eqnarray}\label{eq...}
u(t) \le a G_{\beta}([b\Gamma (\beta)]^{1/{\beta}}t),
\end{eqnarray}
where $G_{\alpha}(x)$ is a monotone function, while $\Gamma (\beta)$ is a
gamma function.

In our case we have
\begin{eqnarray}\label{eq....}
& & \vert v(t) \vert _E \le C(\varepsilon)
G_{\beta}([b\Gamma (\beta)]^{1/{\beta}}t)\le \\
& & C(\varepsilon)
G_{\beta}([b\Gamma (\beta)]^{1/{\beta}}L):=d(\varepsilon)\to 0
\ (\beta := 1-\alpha _1 \in (0,1]) \nonumber
\end{eqnarray}
as $\varepsilon \to 0 $ uniformly w.r.t. $\omega \in \Omega $,
$x\in B[0,R_0]$ and $ t\in [0,L] $ for every $L>0$. The theorem is proved.
\end{proof}

\section{The global averaging principle for the
non-autonomous Navier-Stokes equations}

Let $\Omega$ be a compact metric space, $(\Omega,\mathbb
R,\sigma)$ be a dynamical system on $\Omega$, $E$ be a Banach space and
$\langle E ,\varphi , (\Omega ,\mathbb R ,\sigma) \rangle $ be a cocycle on
$(\Omega ,\mathbb R ,\sigma)$ with fibre $E$.

A family of nonempty compact sets $\{ \ I_{\omega} \ | \ \omega \in \Omega \}
\ (I_{\omega} \subset E)$ is called a local compact attractor
(local compact forward attractor)
if the followings conditions are fulfilled:
\begin{enumerate}
\item
$$I=\bigcup \{I_{\omega}\ :
 \ \omega \in \Omega \}$$
is compact;
\item
$$
\varphi _{\lambda_0} (t,I^{\lambda_0}_\omega,\omega)=
I^{\lambda_0}_{\sigma (t,\omega)}
$$
for all $t\in \mathbb R_+$ and $\omega \in \Omega $;
\item
there exists $R_0>0$ such that $I \subset B(0,R_0):=\{ x \in E | \
\vert x \vert < R_0\} $ and
$$
\lim \limits _{t \to \infty} \sup \limits _{\omega \in \Omega}
\beta (\varphi (t,B[0,R_0],\omega), I)=0
$$
(respectively $\lim \limits _{t \to \infty} \sup \limits _{\omega
\in \Omega} \beta (\varphi (t,B[0,R_0],\omega), I_{\omega t})=0 $)
\end{enumerate}

\begin{theorem}\label{t6.1}
Let $\Lambda $ be a compact metric space, $E$ be a Banach space
and $\varphi _\lambda \ (\lambda \in \Lambda )$ be a cocycle on
$(\Omega,\mathbb R,\sigma)$ with fibre $E$. Suppose that the
following conditions are fulfilled:
\begin{enumerate}
\item
the cocycle $\phi _{\lambda _0}$ admits a local compact forward attractor,
\item
the following relation takes place
\begin{equation}
m_{L}(\lambda):=\sup \{ \vert \varphi _{\lambda }(t,x,\omega)-
\varphi _{\lambda _0}(t,x,\omega) \vert \ : \ 0\le t \le L,
\omega \in \Omega , \vert x \vert \le R_0\}  \to 0
\end{equation}
as $\lambda \to \lambda _0 $ for every positive number $L ;$
\item
every cocycle $\varphi _\lambda $ is asymptotically compact.
\end{enumerate}

Then the next statements are valid:
\begin{enumerate}
\item[a.]
there exists a positive number $\mu $ such that for all $\lambda \in
B[\lambda _0,\mu]:=\{ \lambda \in \Lambda \ : \ \rho (\lambda ,
\lambda _0)\le \mu \}$ the cocycle $\varphi _{\lambda}$ admits in $B[0,R_0]$ 
a
forward attractor $\{I_\omega \ : \ \omega \in \Omega \}$;
\item[b.]
\begin{eqnarray}
\sup \limits _{\omega \in \Omega } \beta (I_{\omega}^{\lambda},
I_{\omega}^{\lambda _0}) \to 0 \nonumber
\end{eqnarray}
as $\lambda \to \lambda _0 .$
\end{enumerate}
\end{theorem}
\begin{proof} Let $\rho >0$ be an arbitrary small number such that
$B[I^{\lambda _0},\rho] \subset B[0,R_0].$ We choose $L=L(\frac{\rho}{3})$ 
according
to the condition
$$
\varphi _{\lambda _0} (t,B[0,R_0],\omega)\subset B[I_{\omega}^{\lambda _0},
\frac{\rho}{3}]
$$
for all $\omega \in \Omega $ and $t\ge L(\frac{\rho}{3}).$ Now we choose
$\varepsilon _0 =\varepsilon _0 (L) $ so that $m(\lambda)<\frac{\rho}{3}$
for all $\lambda \in B[\lambda _0,\varepsilon _0] .$

Let $t_1:=L,$ then we have $\varphi _{\lambda _0}
(t_1,x,\omega) \in
B[I^{\lambda _0}_{\omega t_1},\frac{\rho}{3}] $
and $\varphi _{\lambda} (t_1,x,\omega)\in
B[I^{\lambda_0}_{\omega t_1} $. We take the point $x_1:=
\varphi _{\lambda} (t_1,x,\omega) $ as the initial point and we consider
$\varphi _{\lambda} (t,x_1,\omega t_1) $ on the segment $[0,L], $
$$
\varphi _{\lambda _0}(t,x_1,\omega t_1) ; \
\varphi _{\lambda }(t,x_1,\omega t_1)=\varphi _{\lambda }(t,
\varphi _{\lambda }(t_1,x,\omega ),\omega t_1)=
\varphi _{\lambda }(t+t_1,x,\omega ) .
$$
On this segment $\varphi _{\lambda}(t,x_1,\omega t_1) $ and
$\varphi _{\lambda _0}(t,x_1,\omega t_1) $ will diverge by the value less 
than
$\frac{\rho}{3} .$ Since $\varphi _{\lambda _0}(t,x_1,\omega t_1) \in
B[I_{\omega t_1},\frac{\rho}{3}],$ we get
$ \varphi _{\lambda }(2t_1,x,\omega ) \in B[I_{\omega 
2t_1},\frac{2\rho}{3}].$

If we take the point $ x_2:= \varphi _{\lambda }(2t_1,x,\omega ) $
as the initial one, then we see that the situation is similar to
that occurred at the previous step.

Repeating this process, we arrive at a conclusion that $ \varphi
_{\lambda }(t,x,\omega ) \in B[I^{\lambda _0}_{\omega
t},\rho]\subset B(0,R_0)$ for all $t\ge L(\frac{\rho}{3}) $ and
$\omega \in \Omega .$ Since the cocycle $\varphi _{\lambda}$ is
asymptotical compact then according to Theorem \ref{t3.3} and
corollary \ref{cor3.5} it admits a forward attractor $\{
I_{\lambda} \ : \ \omega \in \Omega \}$ such that $I^{\lambda}:=
\bigcup \{I_{\omega}^{\lambda} \ : \ \omega \in \Omega \}
\subseteq B[I^{\lambda _0}, \rho]$ and, consequently, $\beta
(I^{\lambda}, I^{\lambda _0})\to 0$ as $\lambda \to \lambda _0 .$

Below we proved the inclusion $\varphi _{\lambda}(t,B[0,R_0],\omega) 
\subseteq
B[I^{\lambda_0}_{\omega t},\frac{\rho}{3}]$ for all $t\ge L$ and
$\omega \in \Omega $ and, consequently, we obtain
\begin{equation}\label{eq6.1}
\varphi _{\lambda}(t,B[0,R_0],\omega _t) \subseteq
B[I^{\lambda_0}_{\omega },\frac{\rho}{3}]
\end{equation}
for all $t\ge L$ and $\omega \in \Omega .$ Taking onto consideration that
\begin{equation}\label{eq6.2}
I^{\lambda}_{\omega}= \bigcap \limits _{t\ge 0}
\overline{\bigcup \limits _{\tau \ge t} \varphi _{\lambda}(\tau ,B[0,R_0],
\sigma (-\tau,\omega))}
\end{equation}
from (\ref{eq6.1}) and (\ref{eq6.2}) it follows that $I^{\lambda}_{\omega} 
\subseteq
B[I^{\lambda _0}_{\omega},\rho]$ for all $\omega \in \Omega $ and $\lambda 
\in
B[\lambda _0,\varepsilon_0]$ and, consequently,
$\sup \limits _{\omega \in \Omega} \beta 
(I^{\lambda}_{\omega},I^{\lambda_0}_{\omega})\to 0$
as $\lambda \to \lambda _0 .$ The theorem is proved.
\end{proof}

\begin{remark}\label{r6.2} 1. The second condition of Theorem \ref{t6.1} is 
fulfilled,
for example, if the space $E$ is finite-dimensional and the mapping
$\varphi : \mathbb R_{+} \times E \times \Omega \times \Lambda \to E ,$ 
defined by
the equality $\varphi (t,x,\omega,\lambda ):=\varphi_{\lambda} (t,x,\omega) 
$,
is continuous.

In fact, if we suppose that it is not true, then there exist $L_0>0, \ 
\lambda _k \to \lambda _0, \
x_k \in B[0,R_0], t_l \in [0,L_0]$ and $ \omega _k \in \Omega $ such that
\begin{equation}\label{eq6.3}
\vert \varphi_{\lambda _k} (t_k,x_k,\omega _k) -\varphi_{\lambda_0}
(t_k,x_k,\omega _k)\vert \ge \varepsilon _0 >0 .
\end{equation}
Since the sets $B[0,R_0], \Omega$ and $[0,L_0]$ are compacts, we can suppose 
that
the sequences $ \{x_k\},\ \{t_k\} $ and $\{\omega _k\}$ are convergent. 
Denote
by $t_0:= \lim \limits _{k\to \infty} t_k, \ x_0:= \lim \limits
_{k\to \infty} x_k$ and $\omega _0:= \lim \limits _{k\to \infty} \omega _k .$
Passing to limit in the equality (\ref{eq6.3}) and taking into account the 
continuity
of the mapping $\varphi $ we obtain $0\ge \varepsilon _0 .$ The obtained 
contradiction prove our statement.

2. Under the conditions of Theorem \ref{t6.1} if we suppose that the cocyle
$\varphi _{\lambda _0}$ admits a compact global forward attractor
$ \{I^{\lambda _0}_{\omega} \ : \ \omega \in \Omega \}$, i.e.
$$
\lim \limits _{t\to \infty} \sup \limits _{\omega \in \Omega} \beta (\varphi 
_{\lambda _0}
(t,B[0,R],\omega), I_{\omega t})=0
$$
for every $R>0$, then should be naturally to hope that for the
$\lambda $ sufficiently close to $\lambda _0$ the cocycle $\varphi
_{\lambda}$ also will admits a compact global forward attractor
$\{I^{\lambda}_{\omega} \ : \ \omega \in \Omega \}$ in the small
neighborhood of $I^{\lambda_0}$. Unfortunately, generally
speaking, this assertion is not true.

In fact, let $\varphi _0$ be a cocycle (dynamical system)
generated by the equation $x'=-x$ and $\varphi _\lambda $ be a
cocycle generated by the equation $ x'=-x+\lambda x^{3} \ (\lambda
>0) $. It is clear that the cocycle $\varphi _0 \ (\varphi
_\lambda )$ admits a compact global attractor $I^{0}=\{0\} \
(I^{\lambda}=[-\lambda ^{-1/2},\lambda ^{-1/2}])$. In the small
neighborhood of the attractor $I^0=\{0\}$ the cocycle $\varphi
_\lambda \ ($ for small $\lambda$) admits a local (but not global)
attractor $I^{\lambda}=\{0\}.$
\end{remark}

\begin{theorem}\label{t6.3} Let $\Lambda $ be a compact metric space,
$(\Omega , \mathbb R, \sigma )$ be a dynamical system on the compact metric 
space
$\Omega $, $E$ be a Banach space and $\varphi _{\lambda} \ (\lambda \in 
\Lambda )$
be a cocycle on $(\Omega , \mathbb R, \sigma )$ with fibre $E$. Suppose that 
the
following conditions are fulfilled:
\begin{enumerate}
\item
the cocycle $\varphi _{\lambda _0}$ admits a compact global forward 
attractor;
\item
the following relation takes place
$$
m_{L}(\lambda):=\sup \{\vert \varphi _{\lambda }(t,x,\omega)-
\varphi _{\lambda _0}(t,x,\omega) \vert \ :
\ 0\le t \le L,\ \omega \in \Omega ,\ \vert x \vert \le R_0 \} \to 0
$$
as $\lambda \to \lambda _0 $ for every positive number $L ;$
\item
every cocycle $\varphi _\lambda $ admits a compact global attractor
$\{I^{\lambda}_\omega \ : \ \omega \in \Omega \} ;$
\item
the set $I:=\bigcup \{ I^{\lambda} \ : \ \lambda \in \Lambda \} $ is bounded 
in $E$.
\end{enumerate}

Then the following equality
\begin{eqnarray}
\lim \limits _{\lambda \to \lambda _0} \sup \limits _{\omega \in \Omega }
\beta (I_\omega ^{\lambda}, I_\omega ^{\lambda _0})=0 \nonumber
\end{eqnarray}
is fulfilled and, in particularly,
\begin{eqnarray}
\lim \limits _{\lambda \to \lambda _0}
\beta (I^{\lambda}, I^{\lambda _0})=0 . \nonumber
\end{eqnarray}
\end{theorem}
\begin{proof} Suppose that the conditions of the theorem are fulfilled. 
According to the condition
(iv) there exists a positive number $R_0$ such that $I\subset
B(0,R_0).$ Reasoning as in Theorem \ref{t6.1} for all $\rho >0$ we
will find a $L=L(\frac{\rho}{3})>0$ and $\delta _0=\delta
_0(\rho)>0$ such that
$$
\varphi _{\lambda} (t,I^{\lambda}_{\omega}, \omega) \subseteq
B[I^{\lambda _0}_{\omega t},\rho]
$$
for all $t\ge L$ and $\omega \in \Omega $ and, consequently,
$$
I^{\lambda}_{\omega}=\varphi _{\lambda}(t,I_{\sigma (-t,\omega)},
\sigma (-t,\omega)) \subseteq B[I^{\lambda _0}_{\omega},\rho ]
$$
for all $\omega \in \Omega $ and $\rho (\lambda ,\lambda _0) < \delta _0.$
The theorem is proved.
\end{proof}

\begin{lemma}\label{l3.3*}(\cite{Ch20})
Let $ \Lambda $ be a compact metric space and
$ \varphi : \mathbb T _{+} \times W \times \Lambda \times \Omega \mapsto W $
verifies the following conditions :
\begin{enumerate}
\item
$ \varphi $ is continuous ;
\item
for every $ \lambda \in \Lambda $ the function $ \varphi _{\lambda} =
\varphi (\cdot ,\cdot, \lambda, \cdot ) : \mathbb T_{+} \times W \times
\Omega \mapsto W $  is a continuous cocycle on $ \Omega$ with the
fibre $W$;
\item
the cocycle $ \varphi _{\lambda}$ admits a pullback attractor
$ \{I_{\omega}^{\lambda} \ \vert \ \omega \in \Omega \} $ for every
$ \lambda \in \Lambda $;
\item
the set $ \bigcup \{ I^{\lambda} \ \vert \ \lambda \in \Lambda \} $ is
precompact, where $ I^{\lambda} = \bigcup \{ I_{\omega}^{\lambda} \ \vert
\ \omega \in \Omega \},$
\end{enumerate}
then  the following equality
\begin{equation}\label{eq3.2*}
\lim \limits _{ \lambda \to \lambda _{0} , \omega \to \omega _{0} }
\beta ( I_{\omega}^{\lambda}, I_{\omega _{0}}^{\lambda _{0}} )= 0
\end{equation}
takes place for every $ \lambda _{0} \in \Lambda $ and $\omega _{0}\in \Omega 
$ and
\begin{equation}\label{eq3.3*}
\lim \limits _{\lambda \to \lambda _{0} } \beta (I_{\lambda},
I_{\lambda _{0}})=0
\end{equation}
for every $ \lambda _{0} \in \Lambda $.
\end{lemma}

\begin{lemma}\label{l3.5*} (\cite{Ch20})
Let the conditions of Lemma \ref{l3.3*} and additionally the
following condition be fulfilled:
\begin{enumerate}
\item[5.]
for certain $ \lambda _{0} \in \Lambda  $ the application $ F : \Omega 
\mapsto C(W) $,
defined by equality $ F(\omega)= I^{\lambda _{0}}_{\omega} $
is continuous, i.e.
$ \alpha ( F(\omega),F(\omega _{0})) \to 0 $ if $ \omega \to \omega _{0}  $
for every $ \omega _{0} \in \Omega $, where $ \alpha $ is the full metric of
Hausdorff, i.e. $ \alpha (A,B) = \max \{ \beta (A,B), \beta (B,A) \} .$
\end{enumerate}
Then the equality
\begin{equation}\label{eq3.5*}
\lim \limits _{\lambda \to \lambda _{0} } \sup \limits _{\omega \in \Omega }
\beta ( I^{\lambda}_{\omega}, I^{\lambda _{0}}_{\omega} ) = 0
\end{equation}
takes place.
\end{lemma}

\begin{theorem}\label{t4.5*}(\cite{Ch20})
Let $ W $ possess the property $(S)$ and let the cocycle $ \varphi
$ admit a compact pullback attractor $ \{ I_{\omega}\ \vert \
\omega \in \Omega \} $, then :
\begin{enumerate}
\item
the set $ I_{\omega} $ is connected for every $ \omega \in \Omega $;
\item
if the space $ \Omega $ is connected, then the set $ I = \bigcup \{ 
I_{\omega} \vert
\omega \in \Omega \} $ also is connected.
\end{enumerate}
\end{theorem}

\begin{theorem}\label{t6.4}
Let $\varepsilon \in (0,\varepsilon_0)$, $\Omega$ be compact
and connected and
$\varphi_{\varepsilon} \ $( $\bar{\varphi}_{\varepsilon})$ be a cocycle
generated by the equation (\ref{eq5.1}) (respectively by the equation
(\ref{eq5.16})).

Suppose that the following conditions are fulfilled:
\begin{enumerate}
\item
$B(\omega):=B_0(\omega)+B_1(\omega) \ (\omega \in \Omega),
\ B_0,B_1 \in C(\Omega ,L^{2}(E,F))$;
\item
the bilinear forms $B$ and $B_{0}$ satisfy the condition (\ref{eq2.7});
\item
the average of $B_1(\omega)$ is equal to $0$, i.e.
$\lim \limits _{t\to \infty} \frac{1}{t}\int _0^{t}B_1(\omega s)ds =0$ 
uniformly
w.r.t. $\omega \in \Omega $;
\item
the bilinear form $B_0$ satisfies the condition (\ref{eq2.17});
\item
the cocycles $\varphi _{\varepsilon} \ $ and $\bar{\varphi}_{\varepsilon}$
($\varepsilon \in (0,\varepsilon _0])$) are asymptotically compact .
\end{enumerate}

Then the following statements are true:
\begin{enumerate}
\item[a.]
for every $\varepsilon \in (0,\varepsilon _0])$ and $\omega \in \Omega $
the set
$
I^{\varepsilon}_{\omega}:=\{ x \in E \ : \ \text{the solution}
\\ \varphi _{\varepsilon}(t,x,\omega) \ \text{of equation (\ref{eq5.13}) is 
defined
and bounded on}\ \mathbb R \}
$
(respectively \\ $\bar{I}^{\varepsilon}_{\omega}:=\{ x \in E \ : \
\text{the solution} \ \bar{\varphi}_{\varepsilon}(t,x,\omega) \
\text{of equation (\ref{eq5.16}) is defined and bounded on}\
\mathbb R \} $) is nonempty, compact and connected;
\item[b.]
the cocycle $\varphi_{\varepsilon} \ (\bar{\varphi}_{\varepsilon})$ admits a 
compact
global attractor $\{ I^{\varepsilon}_{\omega} \ : \ \omega \in \Omega \}$
(respectively $\{ \bar{I}^{\varepsilon}_{\omega} \ : \ \omega \in \Omega 
\}$);
\item[c.] the set $I^{\varepsilon}$ (respectively $\bar{I}^{\varepsilon}$) is
compact and connected;
\item[d.]
the set $I:=\bigcup \{ I^{\varepsilon} \ : \ \varepsilon \in 
[0,\varepsilon_0] \}$
(respectively $\bar{I}:=\bigcup \{\bar{I}^{\varepsilon} \ : \
\varepsilon \in [0,\varepsilon_0] \}$, where $\bar{I}^{\varepsilon}:
=\bigcup \{\bar{I}^{\varepsilon}_{\omega} \ : \
\omega \in \Omega \}$)
is compact, where $I^{\varepsilon}:=\bigcup
 \{I^{\varepsilon}_{\omega} \ : \
\omega \in \Omega \}$, $I^{0}=\bar{I}^{0}:=\bigcup \{\bar{I}^{0}_{\omega} \ : 
\ \omega
\in \Omega \}$ and $\{\bar{I}^{0}_{\omega} \ : \ \omega
\in \Omega \}$ is a compact global attractor of equation (\ref{eq5.16}), when 
$\varepsilon =1$;
\item[e.]
$\lim \limits _{\varepsilon \to 0} \beta (I^{\varepsilon},\bar{I}^{0})=0$, 
where
$\beta$ is a semi-distance of Hausdorff;
\item[f.] If a dynamical system $(\Omega ,\mathbb R ,\sigma)$ is periodic, 
i.e. there exists
$\omega_{0} \in \Omega $ such that $\omega _{0}\tau = \omega _{0}$ and 
$\Omega =
\{ \omega _{0}t \ : \ t\in [0,\tau) \}$, then
$$
\lim \limits_{\varepsilon \to 0} \sup \{ \beta (I_{\omega}^{\varepsilon},
I_{\omega}^{0} )\}=0.
$$
\end{enumerate}
\end{theorem}
\begin{proof}
Let $\varepsilon \in (0,\varepsilon_0)$, $\Omega$ be compact and
connected and $\varphi _{\varepsilon} \ $(
$\bar{\varphi}_{\varepsilon})$ be a cocycle generated by the
equation (\ref{eq5.1}) (respectively by the equation
(\ref{eq5.16})), then we have
\begin{eqnarray}\label{eq6*}
\varphi _{\varepsilon}(t,x,\omega)=
\varphi (\varepsilon t,x,\omega ,\varepsilon ) (\text{respectively}\
\bar{\varphi}_{\varepsilon}(t,x,\omega)=
\bar{\varphi} (\varepsilon t,x,\omega ,\varepsilon )         ,,
\end{eqnarray}
for all $t\in \mathbb R_{+}, x \in E$ and $\omega \in \Omega,$
where $\varphi (\cdot,\cdot,\cdot,\varepsilon)$ (respectively
$\bar{\varphi} (\cdot,\cdot,\cdot,\varepsilon)$) is a cocycle
generated by the equation (\ref{eq5.17}) (respectively
(\ref{eq5.18})). From the equality (\ref{eq6*}) it follows that
$\{ I^{\varepsilon}_{\omega} \ : \ \omega \in \Omega\} $
(respectively $\{ \bar{I}^{\varepsilon}_{\omega} \ : \ \omega \in
\Omega\},$) is a compact global attractor of the equation
(\ref{eq5.17}) (respectively (\ref{eq5.18})). Now to finish the
proof of theorem it is sufficient to apply Theorems \ref{t5.3},
\ref{t6.3}, \ref{t4.5*} and Lemmas \ref{l3.3*}, \ref{l3.5*}. The
theorem is proved.
\end{proof}

{\bf Acknowledgment:} The research described in this publication
was made possible in part by the Award  MM1-3016 of the Moldovan
Research and Development Association (MRDA) and the U.S. Civilian
Research \& Development Foundation for the Independent States of
the Former Soviet Union (CRDF). This work is also partially
supported by the NSF Grants DMS-0209326 and DMS-0112351. This
paper was written while the first author was visiting the Illinois
Institute of Technology (Department of Applied Mathematics) in
spring of 2002. He would like to thank people in that institution
for their very kind hospitality.

\end{document}